\documentclass{amsart}
\usepackage[latin1]{inputenc}
\usepackage{amssymb}
\usepackage{amsmath}
\usepackage{enumerate}
\usepackage{epsfig,color}

\usepackage[all]{xy}

\newtheorem{theorem}{Theorem}[section]

\newtheorem{proposition}[theorem]{Proposition}
\newtheorem{corollary}[theorem]{Corollary}
\newtheorem{condition}[theorem]{Condition}
\newtheorem{definition}[theorem]{Definition}

\newtheorem{rem}[theorem]{Remark}

\numberwithin{equation}{section}
\newcommand{\rk}{\mbox{rank}}

\newcommand{\ra}{\rightarrow}

\newcommand{\C }{ \mathbb{C}}

\newcommand{\Z}{\mathbb{Z}}
\newcommand{\Q}{\mathbb{Q}}

\newcommand{\N}{\mathbb{N}}

\newcommand{\Pic}{{\rm Pic}}

\makeatletter
\def\blfootnote{\xdef\@thefnmark{}\@footnotetext}

\begin{document}

\keywords{K3 surfaces, Double covers, Bidouble covers.\\
The author is partially supported by PRIN 2010--2011 ``Geometria delle variet\`a algebriche" and FIRB 2012 ``Moduli Spaces and their Applications".}

\title[Smooth double covers of K3 surfaces]{Smooth double covers of K3 surfaces}
\author{Alice Garbagnati}
\address{Alice Garbagnati, Dipartimento di Matematica, Universit\`a di Milano,
  via Saldini 50, I-20133 Milano, Italia}
\email{alice.garbagnati@unimi.it}

\begin{abstract}
In this paper we classify the topological invariants of the possible branch loci of a smooth double cover $f:X\ra Y$ of a K3 surface $Y$. We describe some geometric properties of $X$ which depend on the properties of the branch locus. We give explicit examples of surfaces $X$ with Kodaira dimension 1 and 2 obtained as double cover of K3 surfaces and we describe some of them as bidouble cover of rational surfaces. Then, we classify the K3 surfaces which admit smooth double covers $X$ satisfying certain conditions; under these conditions the surface $X$ is of general type, $h^{1,0}(X)=0$ and $h^{2,0}(X)=2$. We discuss the variation of the Hodge structure of $H^2(X,\Z)$ for some of these surfaces $X$.
\end{abstract}

\maketitle

\section{Introduction}
One of the main and most famous example of K3 surfaces is given by the Kummer surfaces, $Km(A)$, which are, by definition, desingularization of the quotient of an Abelian surface, $A$, by an involution. Equivalently, one can say that a Kummer surface is a K3 surface which admits a double cover branched along 16 disjoint rational curves. This double cover is $A$ blown up in 16 points.  In \cite{NikKummer}, Nikulin characterized the Kummer surfaces by some properties of their N\'eron--Severi groups. This allows to describe the moduli space of the Kummer surfaces in terms of the moduli space of $L$-polarized K3 surfaces for a certain lattice $L$. Nikulin also observed, in the same paper, that the double cover of a K3 surface branched along the disjoint union of rational curves is birational to a surface with trivial canonical bundle and can be branched either on 16 or on 8 curves. In the latter case the cover surface is birational to a K3 surface. In \cite{GSprojective}, the moduli space of K3 surfaces which admit a double cover branched along the disjoint union of 8 rational curves is described in terms of moduli space of polarized K3 surfaces.

Since a smooth double cover of a K3 surface branched along the disjoint union of rational curves is so restrictive, it seems natural to ask what are the conditions on the branch locus in the more general hypothesis that there exists a double cover of a K3 surface branched along the disjoint union of smooth curves, not necessarily rational. Once one determines these conditions, two questions naturally arise. The first is about the geometry of the double cover: what are the numerical invariants of the double cover? The second fundamental question is about the existence and the moduli of the K3 surfaces which admit a double cover with a branch locus of chosen type. We will discuss both these questions, see Theorem \ref{prop: classification}, Section \ref{sebsect:  existence} and Theorem \ref{theorem: classification gC>2, h=-1 }. 

The general setting of this paper will be the following: $Y$ will be a K3 surface and $f:X\ra Y$ is a double cover branched along the disjoint union of $n$ curves $C$ and $R_i$, $i=1,\ldots,n-1$, where $C$ is the curve with highest genus among $C$ and $R_i$.  We always assume that $X$ is a smooth (not necessarily minimal) surface and we refer to this writing that $f:X\ra Y$ is a smooth double cover of $Y$. In Theorem \ref{prop: classification} we describe the admissible configurations of curves $C$, $R_i$, we observe that the genus of $C$ determines the Kodaira dimension of $X$ and that if $g(C)\neq 1$, there are very few possibilities for $n$. 

The problem of the existence of K3 surfaces admitting certain double covers and of the description of their moduli is addressed in Sections \ref{sebsect:  existence}, \ref{sec: double cover h=-1} and \ref{sec:  double cover h not -1}. 

If $g(C)=0$, the results in \cite{NikKummer} and \cite{GSprojective} discussed at the beginning of this introduction show that $n$ is either 8 or 16. In both these cases there exists a countable number of equidimensional irreducible components of the moduli space of the K3 surfaces admitting a double cover branched along $n$ disjoint rational curves and the dimension of these components is 11 if $n=8$ and is $3$ if $n=16$.

If $g(C)=1$, there is no an upper bound on $n$ (which is in any case greater than 1) and the choice of the topological properties of the branch locus does not determine a lattice which characterizes the K3 surfaces admitting the smooth double cover considered. However, the existence and the geometric description of the surfaces which are double covers of K3 surfaces branched on the disjoint union of $n$ curves $C$ and $R_i$ such that $g(C)=1$ is quite clear. All the surfaces constructed in this way have Kodaira dimension 1. The general results on these surfaces are given in the proof of Theorem \ref{prop: classification} and in Section \ref{sebsect:  existence}. Explicit examples are constructed in Section \ref{sec: examples double cover g(C)=1}. 

If $g(C)>1$, then $n\leq 17$ and there are restrictive conditions on the curves in the branch locus (they are all rational with the exception of $C$, whose genus is congruent to $n\mod 4$). The choice of specific values for $n$ and $g(C)$ (i.e. of the topological properties of the branch locus), determines a specific lattice, which can be used to polarize the K3 surfaces admitting the required double cover. However, the problem of the existence is not trivial at all, at least in case $n>1$. The main problem is that once one knows that the N\'eron--Severi group of a K3 surface admits certain required classes (which is a lattice theoretic condition), one has to prove that these classes represent smooth irreducible curves on the surface, in order to construct the double cover of the K3 surface branched on these curves. One of the main result of this paper is the complete description of the K3 surfaces which admit a smooth double cover branched on $n$ curves such that $g(C)=n-4\geq 2$. This result is stated in Theorem \ref{theorem: classification gC>2, h=-1 }: in Proposition \ref{prop: possible NS} we describe the N\'eron--Severi groups of the K3 surfaces which admit the required double covers and in Propositions \ref{prop: Y_n with NS=Ln1}, \ref{prop: Y_n with NS=Ln2} and \ref{prop: Y_n with NS=Lni bog index} we prove that it suffices to know that the N\'eron--Severi group of a K3 surface is contained in a finite known list of lattices in order to conclude that the double cover exists. So, the Theorem \ref{theorem: classification gC>2, h=-1 } is the analogous of the results proved in \cite{NikKummer} in \cite{GSprojective} under the different assumptions that $g(C)=0$. One of the main differences between the results in the case $g(C)=0$ and the one obtained for $g(C)=n-4\geq 2$ is that in the latter case there is a finite number of possible N\'eron--Severi groups (and in fact these are 16).

Once one proves the existence of the K3 surfaces $Y_n$ admitting double covers branched along $n-1$ rational curves and along the curve $C$ of genus $g(C)=n-4\geq 2$, and thus the existence of the surfaces $X_n$ which are double covers of $Y_n$, one describes the geometric properties of these surfaces. In Sections \ref{subsec: existence h=-1} we describe projective models of the K3 surfaces $Y_n$ both as subvariety of $\mathbb{P}^{n-4}$ with $n-1$ nodes and as double cover of $\mathbb{P}^2$ branched along a reducible sextic. As a byproduct of the latter description one obtains some curious information on the existence of a certain number of curves which cut out an odd theta characteristic on  smooth plane quintic (resp. smooth plane quartic) and which are moreover tangent to a line (resp. to a conic), see Corollary \ref{cor: existence theta characteristics}. 

In two specific cases ($n=6$ and $n=9$) we also provide an easy an self contained description of the surface $X_n$, see Section \ref{subsec: Y6 and  X6}. It turns out that $X_n$ is a bidouble cover of $\mathbb{P}^2$. This result, which is of course of a certain geometric interest, also has important consequences in the description of the variation of the Hodge structure of $X_n$. 
One of the reasons for which we concentrate ourselves on the double cover branched along $n-1$ rational curves and a curve of genus $g(C)=n-4\geq 2$, is that, under this condition, the surface $X_n$ has $h^{2,0}(X_n)=2$.
This is particulary interesting because the existence of the double cover $f:X_n\ra Y_n$ implies a splitting of the Hodge structure of the middle cohomology of $X_n$ in two sub--Hodge structures of K3 type. One of them is induced by the pull back of the Hodge structure on the middle cohomolgy of $Y_n$, the other is its orthogonal complement. In general, this latter Hodge structure of K3 type is not immediately related to a K3 surface by our geometric construction. On the other hand, in the specific cases of the surfaces $X_n$ which are bidouble cover of $\mathbb{P}^2$, and are  described in Section \ref{subsec: Y6 and  X6}, it is immediate to check that both the sub--Hodge structures of K3 type defined on $H^2(X_n,\Z)$ are induced by Hodge structures of K3 surfaces. 
Similar results (on the splitting of the Hodge structure of $H^2(X,\Z)$ in sub-Hodge structures of K3 type which are geometrically associated to K3 surfaces) are provided in Sections \ref{sec: examples double cover g(C)=1} and \ref{sec:  double cover h not -1}, under the assumption that $g(C)=1$ or $g(C)=n=10$ respectively.  \\

{\it {\bf Acknowledgments}. I would like to thank Bert van Geemen for many enlightening discussions and suggestions during the preparation of this paper. It's a pleasure to thank also Matteo Penegini for interesting discussions.}

\section{Preliminaries}
The aim of this section is to recall some results which will be useful in the following and to fix the notation.
\subsection{K3 surfaces} 
Let $Y$ be a smooth complex projective surface.
\begin{definition}
The surface $Y$ is a K3 surface if the canonical bundle of $Y$ is trivial and $h^{1,0}(Y)=0$.
\end{definition}
We recall that given a $\Z$-module $V$, a Hodge decomposition of weigh 2 on $V$ is a decomposition $V\otimes \C=V^{2,0}\oplus V^{1,1}\oplus V^{0,2}$ such that $\overline{V^{p,q}}=V^{q,p}$. We will say that a weigh 2 Hodge decomposition $V\otimes \C=V^{2,0}\oplus V^{1,1}\oplus V^{0,2}$ is of K3 type of $\dim(V^{2,0})=1$.

Since the canonical bundle of $Y$ is trivial, the Hodge decomposition of $H^2(Y,\Z)$ induced by the complex structure on $Y$ is of K3 type. 

The second cohomology group with integer coefficients of a K3 surface $Y$ equipped with the cup product is a lattice.
The lattice $(H^2(Y,\Z),\cup)$ is the unique even unimodular lattice of signature $(3,19)$. So it is isometric to a lattice which does not depend on $Y$ and is denoted by $\Lambda_{K3}$.

The N\'eron--Severi group $NS(Y):=H^2(Y,\Z)\cap H^{1,1}(Y)$ coincides with the Picard group and is an even sublattice of $\Lambda_{K3}$ which has signature $(1,\rho(Y)-1)$, where $\rho(Y):=\rk (NS(Y))$ is the Picard number of $Y$. The transcendental lattice $T_Y$, which is defined as the orthogonal to the N\'eron--Severi group in $H^2(Y,\Z)$, is an even sublattice of $H^2(Y, \Z)$ with signature $(2,20-\rho)$. The standard Hodge decomposition on $H^2(Y,\Z)$ induces a Hodge decomposition on $T_Y$, which is of K3 type.

Let $D\in NS(Y)$ be a nef divisor on $Y$. In the following we will denote by $\varphi_{|D|}:Y\ra\mathbb{P}^{\frac{1}{2}D^2+1}$ the map induced by the linear system $|D|$. 
 
\subsection{Double covers and bidouble covers}
Here we recall  some known properties of double covers between surfaces (see e.g. \cite[Section 17 Chapter I and Section 22 Chapter V]{BHPV}) and the definition of bidouble covers (see e.g. \cite{C}).

Let $X$ and $W$ be two surfaces and let $f:X\ra W$ be a double cover branched along the union of the curves $B_i$. Then the branch locus is $B:=\cup_i B_i$ and clearly the classes of the curves $B_i$ (denoted again by $B_i$) are classes in $\Pic(W)$. The divisor $B:=\sum_iB_i$ is the branch divisor of the double cover. \begin{proposition}
Let $f:X\ra W$ be a double cover branched along the union of the curves $B_i$. Then there exists a divisor $L\in Pic(W)$ such that $2L\simeq B$. If the surface $W$ is smooth, then $X$ is smooth if and only if the branch locus $B$ is smooth. 
\end{proposition}
In the previous situation we will say that the set of curves $\{B_i\}$ is 2-divisible, i.e. a set of curves $\{B_i\}$ on $W$ is 2-divisible if there exists a class $L$ in the Picard group such that $2L\simeq \sum_i B_i$.

\begin{proposition}
Let $W$ be a surface such that $\{B_i\}$ is a 2-divisible set of curves on $W$. Then there exists a unique surface $X$ and a map $f:X\ra W$ such that $f:X\ra W$ is a double cover branched along $\cup_i B_i$.
\end{proposition}

\begin{definition}
Let $f:X\ra W$ be a finite map, generically $4:1$. It is called a bidouble cover of $W$ if $f:X\ra W$ is a Galois cover, with cover group $(\Z/2\Z)^2$.
\end{definition}

\subsection{Divisible sets of disjoint rational curves on K3 surfaces}\label{subsec: divisible subset of rational curves}
The 2-divisible sets of disjoint rational curves on K3 surfaces are well known and studied, see e.g. \cite{NikKummer}, \cite{NikSympl} and here we recall the main results which will be useful in the following.

Let us consider a set of disjoint rational curves $\mathcal{R}_m:=\{R_j\}$, $j=1,\ldots m$ contained in a K3 surface. If the set $\mathcal{R}_m$ is 2-divisible, then either $m=8$ or $m=16$, \cite{NikKummer}. Clearly the set $\mathcal{R}_m$ with $m\geq 8$ can contain some subsets of curves which are 2-divisible and one can ask if it is possible that there are more than one subset of 8 curves which are two divisible. The answer is that this is possible, but if two subsets of 8 disjoint rational curves are 2-divisible, then their intersection is either empty or consists of exactly 4 curves. 

\begin{definition}
In the following we will denote by:\\ 
$\bullet$ $M_{\Z/2\Z}$ the lattice generated by the classes in $\mathcal{R}_8$ and by the class $\left(\sum_{i=1}^8R_i\right)/2$;\\
$\bullet$ 
$M_{(\Z/2\Z)^2}$ the lattice generated by the classes in $\mathcal{R}_{12}$ and by the classes $\left(\sum_{i=1}^8R_i\right)/2$, $\left(\sum_{i=5}^{12}R_i\right)/2$;\\
$\bullet$ $M_{(\Z/2\Z)^3}$ the lattice generated by the classes in $\mathcal{R}_{14}$ and by the classes $\left(\sum_{i=1}^8R_i\right)/2$, $\left(\sum_{i=5}^{12}R_i\right)/2$, $\left(R_1+R_2+R_5+R_6+R_9+R_{10}+R_{13}+R_{14}\right)/2$;\\
$\bullet$ $M_{(\Z/2\Z)^4}$ the lattice generated by the classes in $\mathcal{R}_{15}$
and by the classes $\left(\sum_{i=1}^8R_i\right)/2$, $\left(\sum_{i=5}^{12}R_i\right)/2$, $\left(R_1+R_2+R_5+R_6+R_9+R_{10}+R_{13}+R_{14}\right)/2$,\\ $\left(R_1+R_3+R_5+R_7+R_9+R_{11}+R_{13}+R_{15}\right)/2$;\\
$\bullet$  $K$ the lattice generated by the classes in $\mathcal{R}_{16}$
and by the classes $\left(\sum_{i=1}^8R_i\right)/2$, $\left(\sum_{i=5}^{12}R_i\right)/2$, $\left(R_1+R_2+R_5+R_6+R_9+R_{10}+R_{13}+R_{14}\right)/2$,\\ $\left(R_1+R_3+R_5+R_7+R_9+R_{11}+R_{13}+R_{15}\right)/2$,$\left(\sum_{i=1}^{16}R_i\right)/2$. \end{definition}

The main results on $\mathcal{R}_m$ are summarized in the following proposition:
\begin{proposition}\label{prop: the setes Rm}
Let $\mathcal{R}_m$ be a set of $m$ disjoint rational curves on a K3 surface $W$. Then $m\leq 16$ and the minimal primitive sublattice of $NS(W)$ which contains all the curves in $\mathcal{R}_m$ is one of the followings: \begin{itemize}
\item $\langle -2\rangle^m$ if $m\leq 11$;
\item $M_{\Z/2\Z}\oplus \langle -2\rangle^{m-8}$ if $8\leq m\leq 12$;
\item $M_{(\Z/2\Z)^2}\oplus \langle -2\rangle^{m-12}$ if $12\leq m\leq 13$;
\item $M_{(\Z/2\Z)^3}$ if $m=14$;
\item $M_{(\Z/2\Z)^4}$ if $m=15$;
\item $K$ if $m=16$.
\end{itemize}
\end{proposition}

\section{Smooth double covers of K3 surfaces: admissible branch loci}

In this section we discuss the topological properties of the branch locus of a smooth double cover $f:X\ra Y$ of a K3 surface $Y$. The main result is Theorem \ref{prop: classification} in which we assume the existence of a double cover of a K3 surface, branched along $n$ curves, we denote by $C$ the curve with highest genus among the branch curves, we prove that $n$ and $g(C)$ have to satisfy certain conditions and we describe several properties of $X$ which depend only on the topology of the branch locus of $f:X\ra Y$; in particular we will prove that the genus of $C$ determines the Kodaira dimension of $X$. Since Theorem \ref{prop: classification} is stated under the assumption of that there exists a certain double cover, it says nothing about the existence of these double covers. The problem of the existence is discussed in Section \ref{sebsect:  existence}.

\subsection{Classification of the possible branch loci and general properties of the double cover}

\begin{theorem}\label{prop: classification}
Let $Y$ be a K3 surface and let $f:X\ra Y$ be a smooth double cover branched along the union of $n$ disjoint smooth curves $C$,
$R_1,\ldots, R_{n-1}$, with $g(C)\geq g(R_i)$ for
$i=1,\ldots, n-1$. To obtain a minimal model $X_{min}$ of $X$ one has to contract all the rational curves in the ramification locus. Denoted by $L$, the divisor $L\in NS(Y)$ such that $C+\sum_i R_i\sim 2L$, one has: \begin{align}\label{eq: invariants of X}\begin{array}{ll} K_X=\pi^*(\mathcal{O}_Y(L)), & \chi(X)=4+\frac{1}{2}L^2, \\
p_g(X)=1+h^0(Y,\mathcal{O}_Y(L)),& q(X)=-2-\frac{1}{2}L^2+h^0(Y,\mathcal{O}_Y(L)),\\
c_1^2(X)=2L^2,&
c_2(X)=48+4L^2.\end{array}\end{align} 
We distinguish the three different cases $g(C)=0$, $g(C)=1$ and $g(C)>1$: 
\begin{itemize}\item If $g(C)=0$, then $n$ is either 16 or 8. If $n=16$, $X_{min}$ is obtained by $X$ contracting 16 $(-1)$-curves and it is an Abelian surface. If $n=8$, $X_{min}$ is obtained by $X$ contracting 8 $(-1)$-curves and it is a K3 surface.
\item If $g(C)=1$, let us denote by $k$ the number of curves of genus 1 in $\{C,R_1,\ldots R_{n-1}\}$. Then: $$\begin{array}{ll}L^2=(k-n)/2;& n-k\equiv 0 \mod 4;\\
(n-k)/4\leq 4;& b:=k+(n-k)/4\mbox{ is even.}
\end{array}$$ 
The map $\varphi_{|C|}:Y\ra\mathbb{P}^1$ is a genus 1 fibration and there exists a $2:1$ cover $A\ra\mathbb{P}^1$ branched in $b$ points such that $X$ is the fiber product $Y\times_{\mathbb{P}^1} A$. The induced map $X\ra A$ is a genus 1 fibration, the Kodaira dimension of $X$ is 1 and$$c_1^2(X_{min})=0,\ \chi(X_{min})=4-(n-k)/4,\ c_2(X_{min})=48+3k-3n.$$
\item If $g(C)>1$, then: 
$$\begin{array}{llll} g(R_i)=0,&i=1,\ldots, n-1;& g(C)=n+4h, &h\in \Z,\ h\geq -3;\\
L^2=2h;& n\leq 17;& \rho(Y)\geq n.\end{array}$$ The surface $X$ is of general type and $$c_1^2(X_{min})=g(C)-1,\  \chi(X_{min})=4+h,\ c_2(X_{min})=48+8h-(n-1).$$
\end{itemize}
\end{theorem}
\proof
A K3 surface is a simply connected surface, so every double cover of $Y$ is branched. Since the surfaces $X$ and $Y$ are smooth, the branch locus can not contain isolated points, and the branch locus consists of the union of a certain number, $n$, of curves. Again the smoothness of $X$ implies that all the curves in the branch locus are smooth and they do not intersect.  

The values of $h^{i,0}$, $c_1^2$, $c_2$ of $X$ are computed by \cite[Chapter V, Section 22]{BHPV}.

In order to construct a minimal model of $X$ one has to identify the $(-1)$-curves on $X$. Each rational curve in the branch locus of $f:X\ra Y$ has self-intersection $-2$ (since it is a rational curve on a K3 surface). Its inverse image on $X$ is a rational curve with self intersection $-1$. Hence each rational curve in the ramification locus is a contractible curve in $X$. Let us denote by $\gamma_X:X\ra X_{min}$ the contraction of all these curves. We denote by $\gamma_Y:Y\ra Y_{min}$ the contractions of the rational curves in the branch locus, so we have a singular surface $Y_{min}$ and a commutative diagram 
$$\xymatrix{ X\ar[r]^f_{2:1}\ar[d]_{\gamma_X}&Y\ar[d]_{\gamma_Y}\\
X_{min}\ar[r]^l_{2:1}&Y_{min}}$$
where $l:X_{min}\ra Y_{min}$ is a $2:1$ map branched on the singular points of $Y_{min}$ and the image under $\gamma_Y$ of the curves in the branch locus of $f:X\ra Y$ which are not rational.

Let us assume that there is a $(-1)$-curve $E$ on $X_{min}$. In this case we can contract $E$ to obtain a smooth surface $X_{min}'$ and we can also contract $l(E)$ on $Y_{min}$ in order to obtain a singular surface $Y_{min}'$ and a commutative diagram:
$$\xymatrix{ X_{min}\ar[r]^l_{2:1}\ar[d]&Y_{min}\ar[d]\\
X_{min}'\ar[r]_{2:1}&Y_{min}'}$$

Now there are two possibilities: either $E\subset X_{min}$ is the image of a $(-1)$-curve on $X$ which does not intersect the curves contracted by $\gamma_X$, or $E$ is the image of a curve which is not a $(-1)$-curve on $X$ but intersects at least one curve contracted by $\gamma_X$. We immediately exclude the first case, indeed by construction $E$ is not a ramification curve of the $2:1$ cover $f:X\ra Y$ and by projection formula we can exclude that there exists a $(-1)$-curve on $X$ which is mapped to a non branch curve on $Y$. So $E$ is the image of a curve meeting some of the curves contracted by $\gamma_X$, and thus $l(E)$ passes through at least one singular point of $Y_{min}$. After the contraction of $l(E)$, there exists at least one singularity of $Y_{min}'$ which is worst than an $A_1$-singularity. On the other hand the double cover $X_{min}'\ra Y_{min}'$ produces a smooth surface $X_{min}'$, which is impossible. 

So $X_{min}$ can not admit $(-1)$-curves and it is the minimal model of $X$.

The case $g(C)=0$, i.e. all the curves in the branch locus are rational, was considered in \cite{NikKummer}, where it is proved that the canonical bundle of the minimal model of $X$ is trivial, so the minimal model of $X$ is either an Abelian surface or a K3 surface. In the first case $n=16$ and in the second $n=8$.

If $g(C)=1$, then $\varphi_{|C|}:Y\ra \mathbb{P}^1$ is a genus 1 fibration. Since $CR_i=0$, all the curves in the branch locus are contained in the fibers of $\varphi_{|C|}:Y\ra\mathbb{P}^1$. Then the $k$ curves of genus 1 in the branch locus are fibers of the fibration  $\varphi_{|C|}:Y\ra \mathbb{P}^1$ and the $n-k$ rational curves in the branch locus are components of reducible fibers. Let us denote by $\{P_i\}$ the set of the points in $\mathbb{P}^1$ which are the images of the curves $\{C, R_1,\ldots R_{n-1}\}$ for $\varphi_{|C|}$. The fibers $F_t:=\varphi_{|C|}^{-1}(t)$ for which $t\not\in\{P_i\}$ are disjoint from the branch locus, thus  $f^{-1}(F_t)$ is isomorphic to two disjoint copies of $F_t$ and the restriction of $f:X\ra Y$ to $f^{-1}(F_t)$ is an unramified double cover. Hence $X$ carries a genus 1 fibration induced by $\varphi_{|C|}:Y\ra \mathbb{P}^1$ by a base change of order 2, i.e. there exists a $2:1$ cover $A\ra\mathbb{P}^1$, such that $X:=Y\times_{\mathbb{P}^1}A$. The branch locus of $X\ra Y$ is contained in the fibers over the branch points of $A\ra \mathbb{P}^1$, and so $\{P_i\}$ is the set of branch points of the double cover $A\ra\mathbb{P}^1$.
Each curve of genus 1 in the branch locus is a fiber, so it is mapped to a branch point $P_i$, $i=1,\ldots, k$. The $n-k$ rational curves in the branch locus are contained in reducible fibers. Let us assume that a certain number of rational curves in the branch locus are mapped to the same point $P_{k+j}$. 
If $F_{P_{k+j}}$ is a reducible fiber which contains two components $C_1$ and $C_2$, such that both $C_1$ and $C_2$ have an odd multiplicity in $F_{P_{k+j}}$, then $C_1$ and $C_2$ are contained in the branch locus of the double cover $X\ra Y$. Since the curves in the branch locus do not intersect each other, we exclude that $C_1\cap C_2\neq\emptyset$. This implies that the fiber $F_{P_{k+j}}$ is not of type $I_m$ and thus it is a reducible unstable fiber, i.e. it is of one of the following types $I_m^*$, $IV^*$, $III^*$, $II^*$. Each fiber of these types contains exactly 4 rational components with odd multiplicity (and they are all disjoint), so the number of unstable fibers which project to the branch points $P_{k+j}$ for $j>0$ is $(n-k)/4$. Thus the number of points $P_i$ which are the branch points of the double cover $A\ra\mathbb{P}^1$ is $b:=k+(n-k)/4$. The maximal number of reducible unstable fibers in a genus 1 fibration on a K3 surface is 4, thus $(n-k)/4\leq 4$. Moreover, since $b$ is the number of branch points of a double cover of $\mathbb{P}^1$, $b$ is even. The invariants of $X_{min}$ are computed by \eqref{eq:  invariants of X} and by the fact that the number of the rational curves in the branch locus (and so of $(-1)$-curves on $X$) is $n-k$.  
Let us now prove that the Kodaira dimension of $X$ is 1. Since $X$ is a double cover of a K3 surface, the Kodaira dimension of $X$ is not $-\infty$ and since $X$ admits a genus 1 fibration, the Kodaira dimension of $X$ can not be 2. It remains to exclude that the Kodaira dimension of $X$ is 0. All the $(-1)$-curves in $X$ are contained in fibers, so the fibration $X\ra A$ induces a genus 1 fibration $X_{min}\ra A$. If $X$ would be of Kodaira dimension 0, then $c_2(X_{min})\in\{0,12,24\}$. If $c_2(X_{min})=0$, then $48+3k-3n=0$ and so $n-k=16$. This implies that $(n-k)/4=4$. Since $k\geq 1$ and $b$ is even it follows that $b=k+(n-k)/4\geq 6$. This would imply that the genus of $A$ is at least 2, which is not possible for a genus 1 fibration of Kodaira dimension 0. If $c_2(X_{min})=24$, then $X_{min}$ is a K3 surface and $n-k=8$. In this case $b\geq 4$ and thus $g(A)\geq 1$, which is impossible, since a K3 surface does not admit fibrations with a non rational basis. If $c_2(X_{min})=12$, then $X_{min}$ is an Enriques and $n-k=12$. In this case $b\geq 4$ and thus $g(A)\geq 1$, which is impossible, since an Enriques surface does not admit a fibration with non rational basis.

If $g(C)>1$, then $C^2>0$ and by the Hodge index theorem $R_i^2<0$. Since the curves $R_i$ are smooth and irreducible, $R_i^2=-2$ and $g(R_i)=0$. The number of disjoint rational curves on a K3 surface is at most 16 (cf. \cite{NikKummer}) and it follows that $n-1\leq 16$, i.e. $n\leq 17$. 
The classes of the curves $C$ and $R_i$, $i=1,\ldots, n-1$ are all orthogonal and none of them has a trivial self intersection, then the lattice spanned by these classes has rank $n$ and signature $(1,n-1)$. Since this lattice is embedded in the N\'eron--Severi group of $Y$, it follows that $\rho(Y)\geq n$. We recall that also the class $L:=(C+\sum_{i=1}^{n-1}R_i)/2$ is contained in $NS(Y)$. In particular $L^2\in 2\Z$. But $L^2=\left(2g(C)-2-2(n-1)\right)/4$ and then $g(C)\equiv n\mod 4$, i.e. $g(C)=n+4h$ for a certain integer $h$.
The condition $g(C)>1$, implies $n>1-4h$. On the other hand the condition $n\leq 17$ implies $1-4h< 16$, i.e. $h> -4$.
The computation of $L^2$ depends now on the fact that $C^2=2g(C)-2=2n+8h-2$.  The invariants of $X_{min}$ are computed by \eqref{eq:  invariants of X} and by the fact that the number of the rational curves in the branch locus (and so of $(-1)$-curves on $X$) is $n-1$. In particular $c_1^2(X_{min})>0$ and thus $X_{min}$ is either of general type or a rational surface. But $X$ is a double cover of a K3 surface (and thus in particular $h^{2,0}(X)\geq 1$), so $X$ can not be a rational surface.   
\endproof

\subsection{On the existence of the double covers}\label{sebsect: existence}

An easy and well known example of double covers of K3 surfaces branched along a unique curve $C$ can be constructed as follows: let $Y$ be a K3 surface such that $NS(Y)=\Z D$ for a certain divisor $D$ with $D^2=2h>0$. The divisor $D$ can be chosen to be ample (either $D$ or $-D$ is ample). We set $C\simeq 2D$, so $C^2=8h$. The linear system $|C|$ contains a smooth irreducible curve $C$, whose genus is $g(C)=4h+1$. Obviously $C\simeq 2D$ is a 2-divisible divisor, so there exists a double cover of $Y$ branched along the smooth curve $C$. This proves the existence of the double covers described in  Theorem \ref{prop: classification}, case $g(C)>1$, with $n=1$ and $h>0$. Under these conditions the surface $X$ is minimal and $h^{1,0}(X)=3h$, $h^{2,0}(X)=3+4h$ (because $L=D$). Thus, from now on we will always assume that $n\geq 2$.

The existence of K3 surfaces which admit smooth double covers branched along the union of disjoint rational curves is well known, and essentially due to Nikulin, in \cite{NikKummer} and \cite{NikSympl}.

The existence of K3 surfaces which admit a double cover branched along the disjoint union of a certain number of curves of genus 1 and (possibly) of rational curves, follows immediately by the existence of genus 1 fibrations (possibly with unstable reducible fibers) on certain K3 surfaces. Indeed once one has a genus 1 fibration $\varphi:Y\ra\mathbb{P}^1$ it suffices to consider a double cover of $\mathbb{P}^1$ branched on points $P_i$ such that the fibers of $\varphi$ over $P_i$ are either smooth or unstable reducible fibers. This gives examples of the required double cover. More explicit examples will be considered in the Section \ref{sec: examples double cover g(C)=1}.

The existence of K3 surfaces which admit a smooth double cover branched on the disjoint union of $n\geq 2$ curves such that one of them, denoted by $C$, has genus greater than 1 is less evident. A detailed construction of these surfaces and the proof of their existence under some conditions on $n$ and $h$ will be considered in Sections \ref{sec: double cover h=-1} and \ref{sec:  double cover h not -1}.

\section{Double covers branched on at least one curve of genus 1}\label{sec: examples double cover g(C)=1}
This section is devoted to the K3 surfaces which admit a double cover branched on curves of genus 1 and possibly curves of genus 0. In particular, we will produces examples of fibrations $X\ra A$ for which the values of $(h^{2,0}(X),g(A))$ are $(3,0)$ (Sections \ref{subsec:  NS(Y)=U} and \ref{subsec: NS(Y)=U(2)}), $(2,0)$ (Sections \ref{subsec: NS(Y)=U+D4} and \ref{subsect: NS(Y)=U(2)+D4}) and $(2,1)$ (Sections \ref{subsec: NS(Y)=U+D4+D4} and \ref{subsec: NS(Y)=U(2)+D4+D4})). The main difference between the fibrations considered in Sections \ref{subsec:  NS(Y)=U}, \ref{subsec: NS(Y)=U+D4}, \ref{subsec: NS(Y)=U+D4+D4} and the ones considered in Sections \ref{subsec: NS(Y)=U(2)}, \ref{subsect: NS(Y)=U(2)+D4}, 
\ref{subsec: NS(Y)=U(2)+D4+D4} is the presence or not of a section for the fibration.
We will discuss the possibility to obtain the double cover of the K3 surface as bidouble covers either of $\mathbb{F}_4$ (if the fibration has a section) or of $\mathbb{P}^1\times\mathbb{P}^1$ (if the fibration has a bisection, but not a section).\\

Let us assume that $Y$ is a K3 surface which admits a smooth double cover branched on $n$ disjoint curves $\{ C, R_1,\ldots ,R_{n-1}\}$ such that $k>0$ curves in $\{C, R_1,\ldots ,R_{n-1}\}$ have genus 1 and the remaining $n-k$ are rational curves. The $n-k$ rational curves are all linearly independent and they span a negative definite lattice of rank $n-k$. The lattice spanned by the curves $C,R_1,\ldots R_{n-1}$ has rank $n-k+1$ and it is degenerate (i.e.\ there exists a kernel for the quadratic form associated to this lattice). The reason is that all the genus 1 curves are linearly equivalent, and thus they are represented by a unique class, say $C$, whose self intersection is 0. 

In particular, the lattice spanned by the curves $\{C, R_1,\ldots ,R_{n-1}\}$ can not be the N\'eron--Severi group of a K3 surface. Since $Y$ is a projective K3 surface, there exists a (multi)section for the fibration $\varphi_{|C|}:Y\ra \mathbb{P}^1$. So there exists a divisor $P$ such that $PC=m>0$. While $PR_i=m$ for all $i=1,\ldots, k$, there are several possibilities for the value $PR_i$ for $i>k$. Indeed if $i>k$, $PR_i$ has to be non negative, surely not greater than  $m$, and can depend on $i$. In the following we will make some (arbitrary) choices on $m$ and $k$, in order to describe explicitly some examples which seem to be of particular interest.

\subsection{The K3 surface $Y$ with $NS(Y)\simeq U$}\label{subsec: NS(Y)=U}
We are assuming that on the K3 surface $Y$ there is a smooth irreducible curve of genus 1, denoted by $C$. Then $C^2=0$ and $\varphi_{|C|}:Y\ra\mathbb{P}^1$ is a genus 1 fibration. We now assume that the fibration $\varphi_{|C|}:Y\ra\mathbb{P}^1$ admits a section $P$, so that $P^2=-2$ and $PC=1$.
The lattice spanned by $C$ and $C+P$ is a copy of $U$ and the elliptic fibration on the K3 surface has a Weierstrass form, which is $$y^2=x^3+A_{8}(t:s)x+B_{12}(t:s).$$
For a general choice of $A_8$ and $B_{12}$, homogenous polynomials of degree 8 and 12 respectively, the singular fibers are $24I_1$.
Since this elliptic fibration has no reducible fibers, it is not possible to construct a smooth double cover of $Y$ branched on $C$ and some rational curves, but it is possible to construct a double cover branched along the disjoint union of $k$ smooth genus 1 curve. So, with the notation of the Theorem \ref{prop: classification}, we have $n=k=b$, which is an even number. In this case the surface $X$ is smooth, minimal and it carries an elliptic fibration $X\ra A$, where $A$ is a smooth curve with $g(A)=(k/2)-1$. For general choice of $A_8$, $B_{12}$ the surface $X$ is not a product and $h^{1,0}(X)=g(A)=(k/2)-1$. The divisor $L:=(C+\sum_{i=1}^{k-1}R_i)/2$ is linearly equivalent to $k/2$ times the class of the fiber of the fibration, i.e. $L\simeq (k/2) C$ and $L^2=0$. Then $h^1(Y,L)=k/2$ and so $h^0(Y,L)=1+k/2$. This implies that the birational invariant of $X$ are $h^{1,0}(X)=k/2-1$ and $h^{2,0}(X)=2+k/2$.

The case $k=2$ is of particular interest, since in this case $A\simeq \mathbb{P}^1$ and so the elliptic fibration defined on $X$ has a rational basis. In this case the Weierstrass equation of the elliptic fibration on $X$ can be chosen to be 
\begin{equation}\label{eq: first Weierstrass equation}y^2=x^3+A_{8}(t^2:s^2)x+B_{12}(t^2:s^2).\end{equation}
The singular fibers of this fibration are $48 I_1$. 
The elliptic fibration \eqref{eq: first Weierstrass  equation} clearly admits three commuting involutions (which are in fact involutions of the surface $X$, which preserves the elliptic fibration \eqref{eq: first Weierstrass  equation}):
\begin{eqnarray}\label{eq: involutions weierstrass forms}\begin{array}{llll}\iota_1:&((t:s);(x,y))&\mapsto &((-t:s);(x,y)),\\ \iota_2:&((t:s);(x,y))&\mapsto &((-t:s);(x,-y)),\\ \iota_3:&((t:s);(x,y))&\mapsto &((t:s);(x,-y)).\end{array}\end{eqnarray}
By construction, the quotient $X/\iota_1\simeq Y$.

We observe that $h^{2,0}(X)=3$, and indeed the holomorphic 2-form defined on the elliptic fibration \eqref{eq: first Weierstrass  equation} (in the affine subset $s=1$) are $t^rdt\wedge dx/y$ for $r=0,1,2$. We observe that $\iota_1$ preserves the form $tdt\wedge dx/y$ but not the forms $dt\wedge dx/y$ and $t^2dt\wedge dx/y$ and indeed the quotient $X/\iota_1\simeq Y$ has exactly one holomorphic 2-form, since it is a K3 surface. The involution $\iota_2$ preserves the form $dt\wedge dx/y$ and $t^2dt\wedge dx/y$, but not $tdt\wedge dx/y$, so $X/\iota_2$ is a surface with two 2-holomorphic forms. The third involution, $\iota_3$ does not preserves any 2-holomorphic form, hence $h^{2,0}(X/\iota_3)=0$.

Let us consider the involution $\iota_2$: its fixed locus is contained in the fibers over $t=0$ and $t=\infty$. The fiber over $t=0$ is preserved by $\iota_2$, which acts on this fiber as the elliptic involution, hence on this fiber there are exactly four isolated fixed points. In order to study the action of $\iota_2$ over $t=\infty$, we have to consider the change of coordinates $(t,x,y)\mapsto (T,X,Y):=(1/t, x/t^8, y/t^{12})$. The involution $\iota_2$ acts on $(T,X,Y)$ in the following way $(T,X,Y)\mapsto (-T,X,-Y)$. Now the fiber on $t=\infty$ is the fiber on $T=0$ and also on this fiber $\iota_2$ acts as the elliptic involution and thus it fixes exactly 4 isolated points on this fiber.

We observe that $\iota_2$ switches all the fibers over $t\neq 0$ and $t\neq \infty$ (mapping a point of the fiber over $t$ to its opposite in the fiber over $-t$). Hence $X/\iota_2$ admits an elliptic fibration with rational base with $24I_1$ singular fibers which do not contain singular points of this surface and two fibers such that each of them contains 4 singular points of this surface. The minimal resolution of $X/\iota_2$ is an elliptic fibration with singular fibers $2I_0^*+24I_1$ and thus the Euler characteristic of this surface is 36.
This implies that the Kodaira dimension of $X/\iota_2$ is 1 and $h^{2,0}(X/\iota_2)=2$.

The quotient $X/\iota_3$ carries a natural fibration in rational curves, since the involution $\iota_3$ acts on each fiber of $X\ra \mathbb{P}^1$ as the elliptic involution, and thus the fibration induced on $X/\iota_3$ is a fibration in rational curves, in particular $X/\iota_3$ has Kodaira dimension $-\infty$.

The involution $\iota_3$ (and also $\iota_2$) induces on $Y\simeq X/\iota_1$ the elliptic involution and it is known that the quotient of $Y$ by the elliptic involution is isomorphic to the Hirzebruch surface $\mathbb{F}_4$.    

So $\mathbb{F}_4\simeq X/\langle\iota_1,\iota_3\rangle$, i.e. $X$ is birational to a biduoble cover of $\mathbb{F}_4$ and hence we have the following commutative diagram:
$$\xymatrix{&X\ar[dl]\ar[d]\ar[dr]&\\
X/\iota_1\ar[dr]&X/\iota_2\ar[d]&X/\iota_3\ar[dl]\\
&\mathbb{F}_4&}$$

\subsection{The K3 surface $Y$ with $NS(Y)\simeq U\oplus D_4$}\label{subsec: NS(Y)=U+D4}
Let us assume that $NS(Y)\simeq U\oplus D_4$. In this case $Y$ carries an elliptic fibration with singular fibers $I_0^*+18I_1$.
So it is possible to choose the curve $C$ to be a smooth fiber of this fibration and $R_i$, $i=1,\ldots 4$ as the components of $I_0^*$ which have multiplicity 1. This choice corresponds to $n=5$, $k=1$, $b=2$ and thus it produces a surface $X$ which admits an elliptic fibration over $\mathbb{P}^1$.

Up to a choice of the coordinates of $\mathbb{P}^1$, we can assume that the fiber of type $I_0^*$ is over $t=0$. So, by \cite[Table IV.3.1]{M}, the elliptic fibration on $Y$ has equation
$$y^2=x^3+t^2 A_6(t)x+t^3B_9(t).$$
The Weierstrass form of $X$ is 
\begin{equation}\label{eq: second Weierstrass equation}
y^2=x^3+A_6(t^2)x+B_9(t^2).
\end{equation}
The singular fibers of $X$ are $36 I_1$ and $X$ has the following invariants $h^{1,0}(X)=0$, $h^{2,0}(X)=2$. The 2-holomorphic form of $X$ are $t^rdt\wedge dx/y$, $r=0,1$.

The three commuting involutions \eqref{eq: involutions weierstrass  forms} are defined over $X$. The involution $\iota_1$ preserves the form $tdt\wedge dx/y$, the involution $\iota_2$ preserves the form $dt\wedge dx/y$ and the involution $\iota_3$ does not preserve any holomorphic form on $X$. So both $X/\iota_1$ and $X/\iota_2$ have exactly one 2-holomorphic form induced by the ones of $X$.

Here we consider the action of $\iota_1$ on $X$: its fixed locus is contained in the fiber over $t=0$ and $t=\infty$. On the fiber over $t=0$, the automorphism $\iota_1$ acts as the identity. In order to analyze the action of $\iota_1$ over $t=\infty$, we have to consider the change of coordinates $(t,x,y)\mapsto (T,X,Y):=(1/t, x/t^6, y/t^9)$. The involution $\iota_1$ acts on $(T,X,Y)$ in the following way $(T,X,Y)\mapsto (-T,X,-Y)$. Now the fiber on $t=\infty$ is the fiber on $T=0$ and on this fiber $\iota_1$ acts as the elliptic involution and thus it fixes exactly 4 isolated points on it. 

We observe that $\iota_1$ switches all the fibers over $t\neq 0$ and $t\neq \infty$. Hence $X/\iota_1$ admits an elliptic fibration with rational base, with $18I_1$ singular fibers which do not contain singular points for the surface and one fiber which contains 4 singular points for the surface. The minimal resolution of $X/\iota_1$ is an elliptic fibration with singular fibers $I_0^*+18I_1$.
This implies that the Kodaira dimension of $X/\iota_1$ is 0 and $h^{2,0}(X/\iota_1)=1$, i.e. the minimal resolution of $X/\iota_1$ is a K3 surface.

Similarly one obtains that the minimal model of $X/\iota_2$ carries an elliptic fibrations with a fiber of type $I_0^*$ over $t=0$ and the other singular fibers are $18I_1$, hence the minimal model of $X/\iota_2$ is a K3 surface and it is indeed $Y$.

The K3 surfaces which are the minimal model of $X/\iota_1$ and $X/\iota_2$ carry elliptic fibrations which generically have the same properties (i.e. the Mordell--Weil group is trivial and the singular fibers ar $I_0^*+18 I_1$), hence their N\'eron--Severi groups generically are isometric, and indeed isometric to $U\oplus D_4$.

As in Section \ref{subsec:  NS(Y)=U} the surface $X/\iota_3$ carries a fibration in rational curves and has Kodaira dimension $-\infty$.

This construction is quite similar to the one of the previous section, but the main difference is the value of $h^{2,0}(X)$. In this context $h^{2,0}(X)=2$. The fact that $X$ is a smooth double cover of a K3 surface implies that there is a sub-Hodge structure of $H^2(X,\Z)$ which is of K3 type (the pull-back of the one of $Y$). But since $h^{2,0}(X)=2$ also the orthogonal to this sub-Hodge structure is another sub-Hodge structure of K3 type. Here we proved that also this second Hodge structure of K3 type is the pull back of the Hodge structure of a K3 surface (i.e. of the minimal model $X/\iota_1$).

\subsection{The K3 surface $Y$ with $NS(Y)\simeq U\oplus D_4\oplus D_4$}\label{subsec: NS(Y)=U+D4+D4}

Let us assume that $NS(Y)\simeq U\oplus D_4\oplus D_4$. In this case $Y$ carries an elliptic fibration with singular fibers $2I_0^*+12I_1$.
So it is possible to chose the curve $C$ to be a smooth fiber of this fibration, $R_1$ to be another smooth fiber and $R_i$, $i=2,\ldots 9$ to be the simple components of the fibers $I_0^*$. This choice corresponds to $n=10$, $k=2$, $b=4$ and thus it produces a surface $X$ which admits an elliptic fibration over a genus 1 curve $A$.

The singular fibers of $X$ are $24 I_1$ and $X$ has the following invariants $h^{1,0}(X)=1$, $h^{2,0}(X)=2$. 

The fibration $X\ra A$ has a certain Weierstrass equation of the form $y^2=x^3+C(\tau)x+D(\tau)$, where $C(\tau)$ and $D(\tau)$ are rational functions on $A$. By construction there is an involution $\iota_A:A\ra A$ such that $A/\iota_A\simeq \mathbb{P}^1$. The functions $C(\tau)$ and $D(\tau)$ are invariant under $\iota_A$ and thus there are three commuting involutions defined over $X$:
\begin{eqnarray}\label{eq: involtuions if the base is A}\begin{array}{llll}\iota_1:&(\tau,x,y)&\mapsto& (\iota_A(\tau),x,y)\\ \iota_2:&(\tau,x,y)&\mapsto& (\iota_A(\tau),x,-y)\\ \iota_2:&(\tau,x,y)&\mapsto& (\tau,x,-y).\end{array}\end{eqnarray} These involutions are the analogous of the involutions \eqref{eq: involutions weierstrass  forms} defined in Section \ref{subsec:  NS(Y)=U}. Similarly to section \ref{subsec: NS(Y)=U+D4}, one has that $X/\iota_1$ and $X/\iota_2$ are two surfaces which admit an elliptic fibration over $A/\iota_A\simeq \mathbb{P}^1$ whose singular fibers are $12I_1+2I_0^*$, hence their minimal models are K3 surfaces which generically have N\'eron--Severi group isometric to $U\oplus D_4\oplus D_4$ (and indeed one of them is the K3 surface $Y$). The surface $X/\iota_3$ has a fibration in rational curves with basis $A$, and thus its Kodaira dimension is $-\infty$.

As in Section \ref{subsec: NS(Y)=U+D4} we proved that both the sub--Hodge structures of K3 type of $X$ are obtained by pull back of the Hodge structure of a K3 surface.

\subsection{The K3 surface $Y$ with $NS(Y)\simeq U(2)$}\label{subsec: NS(Y)=U(2)}

Let us assume that $NS(Y)\simeq U(2)$. Then on the K3 surface $Y$ there is a smooth irreducible curve of genus 1, denoted by $C$,  $C^2=0$ and $\varphi_{|C|}:Y\ra\mathbb{P}^1$ is a genus 1 fibration. Moreover the fibration $\varphi_{|C|}:Y\ra\mathbb{P}^1$ does not admit a section, but admits a bisection $P$, so that $PC=2$, which is a smooth curve of genus 1, so that $P^2=0$. Indeed we observe that the lattice spanned by $C$ and $C+P$ is a copy of $U(2)$. The K3 surface $Y$ is a $2:1$ cover of $\mathbb{P}^1\times\mathbb{P}^1$ branched along a curve of bidegree $(4,4)$. Let us denote by $(x_0:x_1)$ the homogeneous coordinates of the first copy of $\mathbb{P}^1$ and by $t$ an affine coordinate on the second copy of $\mathbb{P}^1$. Then one has the following equation for $Y$:
\begin{equation}\label{eq: first double cover of P1xP1}y^2=t^4p_1(x_0:x_1)+t^3p_2(x_0:x_1)+t^2p_3(x_0:x_1)+tp_4(x_0:x_1)+p_5(x_0:x_1)\end{equation}
where $p_i(x_0:x_1)$ are homogenous polynomials of degree 4 in $(x_0:x_1)$. For each fixed value of $t$, the equation \eqref{eq: first double cover of  P1xP1} is the equation of a $2:1$ cover of $\mathbb{P}^1_{(x_0:x_1)}$ branched in 4 points, i.e. the equation \eqref{eq: first double cover of  P1xP1} is the equation of a genus 1 fibration $Y\ra\mathbb{P}^1_t$.

Now we consider the double cover of $Y$ branched along two smooth fibers of the fibration $Y\ra \mathbb{P}^1_t$, e.g. on the fibers over $t=0$ and $t=-\infty$. This double cover is the surface $X$ , which is a minimal surface and admits a genus 1 fibration over $\mathbb{P}^1$ whose equation is  
\begin{equation}\label{eq: first double cover of P1xP1 on X}y^2=t^8p_1(x_0:x_1)+t^6p_2(x_0:x_1)+t^4p_3(x_0:x_1)+t^2p_4(x_0:x_1)+p_5(x_0:x_1)\end{equation}
This fibration has 48 singular fibers, all of type $I_1$, and the invariants of the surface $X$ are $h^{1,0}(X)=0$, $h^{2,0}(X)=3$. 
The involutions \eqref{eq: involutions weierstrass  forms} are well defined also on this genus 1 fibration and we can consider the quotients $X/\iota_1$, $X/\iota_2$ and $X/\iota_3$ as we did in Section \ref{subsec:  NS(Y)=U}. So we obtain that $X/\iota_1\simeq Y$, $X/\iota_2$ is a surface with Kodaira dimension 1 and $h^{2,0}(X/\iota_2)=2$, and $X/\iota_3$ carries a fibration in rational curves, its Kodaira dimension is $-\infty$ and $h^{2,0}(X/\iota_3)=0$.
We obtain the following commutative diagram:
$$\xymatrix{&X\ar[dl]\ar[d]\ar[dr]&\\
X/\iota_1\ar[dr]&X/\iota_2\ar[d]&X/\iota_3\ar[dl]\\
&\mathbb{P}^1\times\mathbb{P}^1&}$$
hence $X$ is a bidouble cover of $\mathbb{P}^1\times\mathbb{P}^1$.
More precisely: the surface $X/\iota_1\simeq Y$ is, by definition, the double cover of $\mathbb{P}^1_{(x_0:x_1)}\times \mathbb{P}^1_{(t:s)}$ branched along the curve of bidegree $(4,4)$ $Z_{4,4}:=V(t^4p_1(x_0:x_1)+t^3sp_2(x_0:x_1)+t^2s^2p_3(x_0:x_1)+ts^3p_4(x_0:x_1)+s^4p_5(x_0:x_1))$.
Denoted by $F_1$ and $F_2$ the curves in $\mathbb{P}^1\times \mathbb{P}^1$ which corresponds to $(t:s)=(0:1)$ and $(t:s)=(1:0)$, the surface $X/\iota_2$ is the $2:1$ cover of $\mathbb{P}^1\times \mathbb{P}^1$ branched along the union of $Z_{4,4}$, $F_1$ and $F_2$ and the surface $X/\iota_3$ is the double cover of $\mathbb{P}^1\times\mathbb{P}^1$ branched along the union of $F_1$ and $F_2$. 

\subsection{The K3 surface $Y$ with $NS(Y)\simeq U(2)\oplus D_4$}\label{subsect: NS(Y)=U(2)+D4}
Let us assume that $NS(Y)\simeq U(2)\oplus D_4$. In this case $Y$ is a double cover of $\mathbb{P}^1\times \mathbb{P}^1$ and the genus 1 fibration induced by the projection on the second factor has a fiber of type $I_0^*$. This corresponds to a specialization of the equation \eqref{eq: first double cover of  P1xP1}: we assume that the branch locus of the double cover $Y\ra\mathbb{P}^1\times\mathbb{P}^1$ splits in the union of a smooth curve of bidegree $(4,3)$ and a smooth curve of bidegree $(0,1)$. So the equation of $Y$ is
\begin{equation}y^2=t\left(t^3p_1(x_0:x_1)+t^2p_2(x_0:x_1)+tp_3(x_0:x_1)+p_4(x_0:x_1)\right)\end{equation}
where $p_i(x_0:x_1)$ are homogeneous polynomials of degree 4 in $(x_0:x_1)$. 
Hence, in analogy with what we did in Section \ref{subsec: NS(Y)=U(2)}, the double cover $X$ has the following equation:
\begin{equation}w^2=t^6p_1(x_0:x_1)+t^4p_2(x_0:x_1)+t^2p_3(x_0:x_1)+p_4(x_0:x_1)\end{equation}
where $w:=yt$.
The surface $X$ carries a genus 1 fibration with singular fibers $36I_1$ and with basis $\mathbb{P}^1_t$. So everything is analogous to the case studied in Section \ref{subsec: NS(Y)=U+D4}, except the fact that here the fibrations have no sections.

\subsection{The K3 surface $Y$ with $NS(Y)\simeq U(2)\oplus D_4\oplus D_4$}\label{subsec: NS(Y)=U(2)+D4+D4}
If $NS(Y)\simeq U(2)\oplus D_4\oplus D_4$, the K3 surface $Y$ is a $2:1$ cover of $\mathbb{P}^1\times\mathbb{P}^1$ branched along the union of two curves of bidegree $(0,1)$ and one curve of bidegree $(4,2)$ and the equation of the K3 surface in this case is 
\begin{equation}\label{eq: double cover of P1xP1,2I0*}y^2=t\left(t^2p_2(x_0:x_1)+tp_3(x_0:x_1)+p_4(x_0:x_1)\right)\end{equation}
where $p_i(x_0:x_1)$ are homogenous polynomials of degree 4 in $(x_0:x_1)$. The fibers over $t=0$ and $t=\infty$ are of type $I_0^*$, the other singular fibers are $12I_1$.
 
Considering, as in Section \ref{subsec: NS(Y)=U+D4+D4}, the double cover of $Y$ branched along two smooth fibers of the fibration and the eight smooth rational curves with multiplicity 1 in the fibers of type $I_0^*$, one obtains a surface $X$ which has a genus 1 fibration whose base is a genus 1 curve $A$. The birational invariants of $X$ are $h^{1,0}(X)=1$ and $h^{2,0}(X)=2$. The involutions \eqref{eq: involtuions if the base is  A} are well defined and the quotient surfaces $X/\iota_1$, $X/\iota_2$ and $X/\iota_3$ have properties analogous to the corresponding surfaces described in Section \ref{subsec: NS(Y)=U+D4+D4}.

\section{The case $g(C)>1$, $h=-1$}\label{sec: double cover h=-1}
In this section we discuss the double covers $f:X_n\ra Y_n$ of K3 surfaces with the following properties (corresponding to the case $g(C)>1$, $h=-1$ of Theorem 
\ref{prop: classification}):
\begin{condition}\label{condition} $Y_n$ is a K3 surface and  $f:X_n\ra Y_n$ is a smooth double cover branched on $n$ smooth disjoint curves such that $n-1$ branch curves, denoted by $R_i$, are rational and the remaining curve, denoted by $C$, has genus $g(C)=n-4\geq 2$. Moreover $\rho(Y_n)=n$.\end{condition}

One of the main results of this section is the existence and the classification of the K3 surfaces $Y_n$ satisfying Condition \ref{condition}. To prove this result we show that the geometric Condition \ref{condition} on $Y_n$ is equivalent to a lattice theoretic condition on the N\'eron--Severi group of $Y_n$, indeed in Theorem \ref{theorem: classification gC>2, h=-1 } we will prove that $Y_n$ is as in Condition \ref{condition} if and only if its N\'eron--Severi group is contained in the finite list of lattice given in Definition \ref{defi: the list of lattices}.  The proof of Theorem \ref{theorem: classification gC>2, h=-1 } is essentially contained in two Sections: in Section \ref{subsec: possible NS} one produces the list of lattices and the methods of that section are lattice theoretic. In Section \ref{subsec: existence h=-1} the K3 surfaces whose N\'eron--Severi group is contained in the list are described from a geometric point of view, and this allows to conclude that they coincide with the surfaces satisfying the Condition \ref{condition}.

Section \ref{subsec: bidouble covers P2} is devoted to two specific examples, in which the cover surface $X_n$ can be described as bidouble cover of $\mathbb{P}^2$. The results of this section are related with the following problem:
if $Y_n$ and $X_n$ are as in Condition \ref{condition}, then $h^{2,0}(X_n)=2$ and there exists an involution, the cover involution, denoted by $\iota_X$, such that $X_n/\iota_X\simeq Y_n$. The involution $\iota_X^*$ acts on $H^2(X_n,\Z)$ and splits the Hodge structure of $H^2(X_n,\Z)$ in two sub-Hodge structures:  one is induced on $H^2(X_n,\Z)^{\iota}$ and the other on $\left(H^2(X_n,\Z)^{\iota}\right)^{\perp}$. Since $H^2(X_n,\Q)^{\iota_X}\simeq H^2(Y_n,\Q)$, the Hodge decomposition on $H^2(X_n,\Q)^{\iota_X}$ is of K3 type, i.e. $H^{(2,0)}(X_n)^{\iota_X}$ has dimension 1. But then also $\left(H^{(2,0)}(X_n)^{\iota_X}\right)^{\perp}$ has dimension 1,  so the Hodge structure on $\left(H^2(X_n,\Z)^{\iota}\right)^{\perp}$ is of K3 type. Clearly the variation of the sub--Hodge structure on  $H^{2}(X_n,\Z)^{\iota_X}$ depends on the variation of the Hodge structure on $H^2(Y_n,\Z)$, but one can ask if also the variation of the sub--Hodge structure on $\left(H^2(X_n,\Z)^{\iota}\right)^{\perp}$ is related to the one of a K3 surface. In Section \ref{subsec: bidouble covers P2} we show that at least for two different surfaces $X_n$ the answer to this question is yes, and indeed the Hodge structure on $\left(H^2(X_n,\Z)^{\iota}\right)^{\perp}$ is induced by the one of a K3 surface (which is in general not isomorphic to $Y_n$).

\subsection{General facts}
\begin{proposition}
If $Y_n$ is a K3 surface as in Condition \ref{condition}, then $n\geq 6$ and the birational invariants of $X_n$ are $\chi(X_n)=3$, $h^{1,0}(X_n)=0$, $h^{2,0}(X_n)=2$. In particular they do not depend on $n$.\end{proposition}
\proof Let $Y_n$ be as in Condition \ref{condition}. The condition $g(C)>1$, implies that $n>5$. The class $L:=\left(C+\sum_{i=1}^{n-1}R_i\right)/2$ has self intersection $L^2=-2$, since $C^2=2n-10$. Thus either $L$ or $-L$ is effective. But $LC=n-5>0$, so $L$ is effective.  For every effective divisor $V$ with self intersection $-2$ on a K3 surface, the space $h^0(V)$ has dimension 1 (indeed if $V$ would be linearly equivalent to an effective divisor $V'$, then $V'$ should be a fixed component of $|V|$, since $VV'=-2$ and thus $V-V'$ should be an effective divisor. But $V-V'=0$). It follows that $h^0(Y,L)=1$ and applying \eqref{eq:  invariants of X} one finds the statement.
\endproof

\begin{proposition}\label{prop: the effective class D}
If $Y_n$ is a K3 surface as in Condition \ref{condition}, then there exists an effective class $D$ with $D^2=-2$ on $Y_n$ such that $C\simeq 2D+\sum_{i=1}^{n-1} R_i$ and $DR_i=1$.

Viceversa if $Y$ is a K3 surface such that $\rho(Y)=n\geq 6$ and there exists on $Y$ $n$ rational curves $D$, $R_i$, $i=1,\ldots n-1$ such that $DR_i=1$, $R_iR_j=0$ if $i\neq j$ and a smooth irreducible curve $C$ such that $C\simeq 2D+\sum_{i=1}^{n-1}R_i$, then $Y$ admits a double cover satisfying the Condition \ref{condition}.
\end{proposition}
\proof
Let $Y_n$ be a surface as in Condition \ref{condition}. Then $L:=(C+\sum_{i=1}^{n-1}R_i)/2$ is a class in $NS(Y_n)$. So also $L-\sum_{i=1}^{n-1}R_i=(C-\sum_{i=1}^{n-1}R_i)/2$ is a class in $NS(Y_n)$ and we call this class $D$. We observe that $D^2=-2$, $DC=n-5>0$ and $DR_i=1$. In particular $D$ is an effective divisor. 

Let $Y$ be a surface which contains the curves $D$, C and $R_i$. Since $C\simeq 2D+\sum_{i=1}^{n-1}R_i$, $C+\sum_{i=1}^{n-1}R_i\simeq 2(D+\sum_{i=1}^{n-1}R_i)$ and thus $\{C, R_i\}$ is a 2-divisible set of disjoint curves. Hence there exists a smooth double cover of $Y$ branched on the union of $C$ and $R_i$. Since $C^2=2n-10$, the smooth curve $C$ has genus $g(C)=n-4$. So the double cover of $Y$ branched on the union of $C$ and the $R_i$'s is as in Condition \ref{condition}. 
\endproof
\begin{rem}{\rm
In Corollary \ref{cor: D is irreducible all cases} we will be able to prove a stronger result on $D$: the effective class $D$ necessarily corresponds to a smooth irreducible rational curve.}\end{rem}

\begin{rem}{\rm
If we assume that on a K3 surface there exist $n$ rational curves $D$ and $R_i$, $i=1,\ldots, n-1$ with $DR_i=1$, $R_iR_j=0$, $i\neq j$ and a smooth irreducible curve $C$ such that $C\simeq 2D+\sum_{i=1}^{n-1}R_i$, then we are essentially saying that in the linear system of $C$ there is a reducible curve which splits in the union of $n$ rational curves $D$ and $R_i$. In Proposition \ref{prop: the effective class D} we assumed that $n\geq 6$, in order to obtain a double cover as in Condition \ref{condition}. In the case $n=5$ the curve $C$ is a smooth curve of genus 1, and defines a genus 1 fibration $\varphi_{|C|}:Y\ra\mathbb{P}^1$. The splitting of a curve in the linear system of $C$ as $2D+\sum_{i=1}^4R_i$ implies that the genus 1 fibration $\varphi_{|C|}:X\ra \mathbb{P}^1$ has a fiber of type $I_0^*$. So, in a certain sense, we are generalizing the construction considered in Sections \ref{subsec: NS(Y)=U+D4}, \ref{subsec: NS(Y)=U+D4+D4}, \ref{subsect: NS(Y)=U(2)+D4} and \ref{subsec: NS(Y)=U(2)+D4+D4}.

In case $n<5$, we have $\left(2D+\sum_{i=1}^{n-1} R_i\right)^2<0$, so there is no smooth irreducible curve $C$ linearly equivalent to $2D+\sum_{i=1}^{n-1} R_i$.
}\end{rem}

\subsection{The N\'eron--Severi group of the surfaces $Y_n$}\label{subsec: possible NS}

In this section we investigate the possibilities for $NS(Y_n)$: our goal will be to construct a list $\mathcal{L}$ of all the lattices such that if $Y_n$ is a K3 surface as in Condition \ref{condition}, then $NS(Y_n)\in \mathcal{L}$.

Let $Y_n$ be as in Condition \ref{condition}. By Proposition \ref{prop: the effective class D} the classes $D$ and $R_i$ are contained in $NS(Y_n)$, where $D^2=R_i^2=-2$, $DR_i=1$, $R_iR_j=0$ if $i\neq j$. So the lattice spanned by $D$ and $R_i$ is contained in $NS(Y_n)$. Since both $NS(Y_n)$ and $\langle D, R_i\rangle$ have rank $n$, the inclusion $\langle D, R_i\rangle\hookrightarrow NS(Y_n)$ has finite index.

Let $L_n^{(1)}$ be the lattice which is represented by the matrix
\begin{eqnarray}\label{eq: NS K3 with double cover}L_n^{(1)}\simeq\left[\begin{array}{rrrrrr} -2&1&1&\ldots &1\\
1&-2&0&\ldots&0\\
1&0&-2&\ldots&0\\
\vdots&\vdots&\vdots&\ddots&\vdots\\
1&0&0&\ldots &-2
\end{array}\right].\end{eqnarray}
It is generated over $\Z$ by the classes $d, r_1,\ldots r_{n-1}$ such that $d^2=r_i^2=-2$, $dr_i=1$, $r_ir_j=0$, $i\neq j$.

Equivalently we can say that $L_n^{(1)}$ is generated over $\Q$ by the vectors $c$, $r_i$, $i=1,\ldots, n-1$ such that $cr_i=r_ir_j=0$ if $i\neq j$, $c^2=2n-10$, $r_i^2=-2$ and in order to have a $\Z$ basis we add the class $l:=(c+\sum_{i=1}^{n-1} r_i)/2$. We observe that $d=l-\sum_{i=1}^{n-1}r_i$ is such that $d^2=-2$ and $dr_i=1$, $i=1,\ldots,n-1$. 

So $L_n^{(1)}$ is a hyperbolic lattice of rank $n$ with $d(L_n^{(1)})=(-1)^{n-1}2^{n-2}(n-5)$ and $(L_n^{(1)})^{\vee}/L_n^{(1)}\simeq \Z/(2n-10)\Z\times (\Z/2\Z)^{n-3}$.

We also define some other lattices, which are overlattices of finite index of $L_n^{(1)}$. 

If $n\geq 9$, then we define $L_n^{(2)}$ as the overlattice of index 2 of $L_n^{(1)}$ obtained adding to a $\Z$ basis of $L_n^{(1)}$ the class $(\sum_{i=1}^{8}r_i)/2$. We have $d(L_n^{(2)})=(-1)^{n-1}2^{n-4}(n-5)$ and if $n\neq 9$, then $(L_n^{(2)})^{\vee}/L_n^{(2)}\simeq \Z/(2n-10)\Z\times (\Z/2\Z)^{n-5}$; if $n=9$, then $(L_9^{(2)})^{\vee}/L_9^{(2)}\simeq (\Z/2\Z)^{7}$.

If $n\geq 13$, then we define $L_n^{(4)}$ as the overlattice of index 2 of $L_n^{(2)}$ obtained adding to a $\Z$ basis of $L_n^{(2)}$ the class $(\sum_{i=5}^{12}r_i)/2$. We observe that $L_n^{(4)}$ is an overlattice of index $2^2$ of $L_n^{(1)}$. We have  $d(L_n^{(4)})=(-1)^{n-1}2^{n-6}(n-5)$ and  $(L_n^{(4)})^{\vee}/L_n^{(4)}\simeq \Z/(2n-10)\Z\times (\Z/2\Z)^{n-7}$.

If $n\geq 15$, then we define $L_n^{(8)}$ as the overlattice of index 2 of $L_n^{(4)}$ obtained adding to a $\Z$ basis of $L_n^{(4)}$ the class $(r_1+r_2+r_5+r_6+r_{9}+r_{10}+r_{13}+r_{14})/2$. We observe that $L_n^{(8)}$ is an overlattice of index $2^3$ of $L_n^{(1)}$. We have  $d(L_n^{(8)})=(-1)^{n-1}2^{n-8}(n-5)$ and  $(L_n^{(8)})^{\vee}/L_n^{(8)}\simeq \Z/(2n-10)\Z\times (\Z/2\Z)^{n-9}$.

If $n\geq 16$, then $L_n^{(16)}$ is the overlattice of index 2 of $L_n^{(8)}$ obtained adding to a $\Z$ basis of $L_n^{(8)}$ the class $(r_1+r_3+r_5+r_7+r_{9}+r_{11}+r_{13}+r_{15})/2$. We observe that $L_n^{(16)}$ is an overlattice of index $2^4$ of $L_n^{(1)}$. We have  $d(L_n^{(16)})=(-1)^{n-1}2^{n-10}(n-5)$ and  $(L_n^{(16)})^{\vee}/L_n^{(16)}\simeq \Z/(2n-10)\Z\times (\Z/2\Z)^{n-11}$.

\begin{proposition}\label{prop: possible NS}
Let $Y_n$ be as in Condition \ref{condition}. Let the lattices $L_n^{(r)}$ be as above.

Then $n\leq 16$ and $NS(Y_n)$ is an overlattice of finite index (possibly 1) of $L_n^{(1)}$. In particular:
\begin{itemize}
\item if $n=6,7,8$, then $NS(Y_n)\simeq L_n^{(1)}$;
\item if $n=9,10,11,12$, then $NS(Y_n)$ is isometric to either $L_n^{(1)}$ or $L_{n}^{(2)}$;
\item if $n=13$, then $NS(Y_n)$ is isometric to either $L_{13}^{(2)}$ or $L_{13}^{(4)}$;
\item if $n=14$, then $NS(Y_n)\simeq L_{14}^{(4)}$;
\item if $n=15$ then $NS(Y_n)\simeq L_{15}^{(8)}$;
\item if $n=16$ then $NS(Y_n)\simeq L_{16}^{(16)}$.
\end{itemize}
\end{proposition}
\proof 

We assume that $Y_n$, $C$ and $R_i$ are as in Condition \ref{condition}, so the classes $c:=[C]$ and $r_i:=[R_i]$ are contained in $NS(Y_n)$, as well as the class $\frac{1}{2}([C]-\sum_i[R_i])=\frac{1}{2}(c-\sum_{i=2}^n r_i)=d$. This implies that the lattice $L_n^{(1)}$ is contained in $NS(Y_n)$. Since this lattice has rank $n$ and we assume that $\rho(Y_n)=n$, this inclusion has a finite index. We now construct the possible overlattices of finite index of $L_n^{(1)}$ which can be admissible N\'eron--Severi groups of $Y_n$.

The discriminant group of $L_n^{(1)}$ is generated by the $n-2$ vectors $\frac{c}{2n-10}+\frac{r_1}{2}$, $\frac{r_i+r_{i+1}}{2}$, $i=1,\ldots, n-3$ and thus an overlattice of $L_n^{(1)}$ is obtained by adding a linear combination $$w=\alpha_1\left(\frac{c}{2n-10}+\frac{r_1}{2}\right)+\left(\sum_{j=1}^{n-3}\alpha_j\left(\frac{r_i+r_{i+1}}{2}\right) \right)$$ such that $w^2\in 2\Z$, to the lattice $L_n^{(1)}$.

So we are assuming that $w\not\in L_n^{(1)}$, but $w\in NS(Y_n)$. Now we distinguish two possibilities: either $2w\in L_n^{(1)}$ or $2w\not\in L_n^{(1)}$. If $2w\in L_n^{(1)}$, then $w$ generates a copy of $\Z/2\Z$ in $NS(Y_n)/L_n^{(1)}$ and thus it is equivalent modulo $L_n^{(1)}$ to $\left(\beta_0 c+\sum_{i=1}^{n-1}\beta_i r_i\right)/2$ with $\beta_i\in \{0,1\}$. If $2w\not\in L_n^{(1)}$, then $2w=(2\alpha_1 c)/(2n-10)$ is a non trivial element in $NS(Y_n)/L_n^{(1)}$. This implies that $c$ is a divisible element in $NS(Y_n)$, i.e. there exists $a\in \N$ such that $c/a\in NS(Y_n)$ and thus $2a^2|c^2=2n-10$. This can happen only if $(n,a)=(9,2)$, $(13,2)$, $(14,3)$, $(17,2)$, so for $n=9,13,14,17$ we have to discuss also this possibility (in the last part of the proof). 

Let us first assume that $c$ is not divisible in $NS(Y_n)$ so that $w$ is of the form $(\beta_0 c+\sum_{j=1}^{n-1} \beta_j r_j)/2, \mod L_n^{(1)}$ where $\beta_i\in \{0,1\}$. 

If $\beta_0\equiv 0\mod 2$, then we are adding a linear combination of the classes of the rational curves $R_i$ with coefficients in $\frac{1}{2}\Z$ and the class $w$ satisfies the conditions on the 2-divisible set of disjoint rational curves given in Proposition \ref{prop: the setes Rm}. If we assume $n\leq 16$, thanks to results recalled in Proposition \ref{prop: the setes Rm}, the unique overlattices of $L_n^{(1)}$ obtained adding classes $w$ of the form $(\sum_{j=1}^{n-1} \beta_j r_j)/2$, $\beta_j\in \Z/2\Z$, are the lattices $L_n^{(2)}$, $L_n^{(4)}$, $L_n^{(8)}$, $L_n^{(16)}$ defined above. 

If $\beta_0\equiv 1\mod 2$, then the class $d-w=(c-\sum_{i=1}^{n-1} r_i)/2-w\mod L_n^{(1)}$ is a linear combination of the classes of the rational curves $R_i$ with coefficients in $\frac{1}{2}\Z$, so $d-w$ is a vector of type $(\sum_{j=1}^{n-1} \beta_j r_j)/2$, $\beta_j\in \Z/2\Z$. But clearly the lattice obtained adding the vector $w$ is the same lattice obtained adding the lattice $d-w$ and so we are exactly in the previous case (i.e. the overlattice is obtained adding vectors of type $(\sum_{j=1}^{n-1} \beta_j r_j)/2$, $\beta_j\in \Z/2\Z$). Hence if we assume $n\leq 16$, we obtain again the lattices  $L_n^{(2)}$, $L_n^{(4)}$, $L_n^{(8)}$, $L_n^{(16)}$ defined above.  

In order to prove that $n\leq 16$ and the other results on the possible N\'eron--Severi groups, we have to recall the following well known property on the sublattices of $\Lambda_{K3}$: let $r(L)$ be the rank of a lattice $L$ which is primitively embedded in $\Lambda_{K3}$ (for example the N\'eron--Severi lattice of a K3 surface) and $l(L)$ be its length, then $l(L)\leq \min\{r(L),22-r(L)\}$.

Let us consider the lattice $L_n^{(1)}$. As we said, its discriminant group is generated by the $n-2$ vectors $\frac{c}{2n-10}+\frac{r_1}{2}$, $\frac{r_i+r_{i+1}}{2}$, $i=1,\ldots, n-3$. So the rank is $r(L_n^{(1)})=n$ and the length is $l(L_n^{(1)})=n-2$. This lattice can be the N\'eron--Severi of a K3 surface only if $n-2\leq \min\{n,22-n\}$. In particular if $n>13$ then $NS(Y_n)$ can not be isometric to $L_n^{(1)}$. 

Let us now assume that $n\geq 9$. In this case we have a set of $n-1$ disjoint rational curves which can contain a (sub)set of 8 disjoint rational curves which is 2-divisible. Let us assume that there is exactly one (sub)set of 8 disjoint rational curves which is 2-divisible and (up to a choice of the indices) we assume that it is $\{r_1,\ldots r_8\}$. So we add to the basis of $L_n^{(1)}$ the class $\sum_{i=1}^8r_i/2$ obtaining the overlattice $L_n^{(2)}$. Here we have to consider two cases: $n=9$ and $n>9$. If $n=9$, then $L_n^{(2)}$ contains the vector $v=(c+\sum_{i=1}^8r_i)/2$  and the vector $\sum_{i=1}^8 r_i/2$, thus it contains also the vector $c/2$ (and this is the case where $c$ is divisible). The basis of the discriminant group of $L_9^{(2)}$ is given by $(c/2)/2$, $(r_i+r_{i+1})/2$, $i=1,2,3,4,5,6$. The length is $l(L_9^{(2)})=7$.
If $n>9$, then the discriminant group of $L_9^{(2)}$ is generated by $c/(2n-10)+r_{9}/2$, $(r_i+r_{i+1})/2$, $i\leq n-3$, $i\neq 8,9$. In this cases the length is $l(L_n^{(2)})=n-4$.  As in case $L_n^{(1)}$, also in this case the comparison between $r(L_n^{(2)})=n$ and $l(L_n^{(2)})$ gives condition on $n$. If $n>14$, then $NS(Y_n)$ can not be isometric to $L_n^{(2)}$.

Similarly if $n\geq 13$ we have a set of $n-1$ disjoint rational curves which can contain two subsets of 8 disjoint rational curves which are 2-divisible and up to a permutation of the indices this can be done in a unique way. So let us add to the basis of $L_n^{(2)}$ also the class $\sum_{i=5}^{12}r_i/2$. We obtain the lattices $L_n^{(4)}$. The length of this lattice can be computed as before and it is $l(L_n^{(4)})=n-6$. So, if $n>15$, then $NS(Y_n)$ can not be isometric to $L_n^{(4)}$.

If $n\geq 15$  (resp. 16) we have a set of $n-1$ disjoint rational curves which can contains 3 (resp. 4) subsets of 8 disjoint rational curves which are 2-divisible and up to a permutation of the indices this can be done in a unique way. So let us add the the basis of $L_n^{(4)}$ (resp. $L_n^{(8)}$) another 2-divisible class. We obtain the lattice $L_n^{(8)}$ (resp. $L_n^{(16)}$). The length of this lattice can be computed as before and it is $l(L_n^{(8)})=n-8$ (resp. $l(L_n^{(16)})=n-10$). So, if $n>16$ (resp. $n>17$),  then $NS(Y_n)$ can not be isometric to $L_n^{(8)}$ (resp. $L_n^{(16)}$).

If $n= 17$, the K3 surface $Y_n$ should contains 16 disjoint rational curves. In this case $Y_{17}$ should be a Kummer surface and $NS(Y_{17})$ contains also the sum of all these classes divided by 2. But then it contains also the class $c/2$ ($=(c-\sum_{i=1}^{16}r_i)/2+\sum_{i=1}^{16}t_i/2)$. Its self-intersection should be $c^2/4=(34-10)/4=6$. Since $c/2$ is orthogonal to the classes $r_j$, $NS(Y_{17})$ should be an overlattice of finite index of $\langle 6\rangle\oplus K$ in which $K$ is primitively embedded. But there are no admissible overlattice of $\langle 6\rangle\oplus K$ in which $K$ is primitively embedded and which has index greater than 1, so $NS(Y_{17})$ should be equal to $\langle 6\rangle\oplus K$, which is impossible, because the length of $\langle 6\rangle\oplus K$ is 7 and its rank is 17.

If $n\geq 14$ (resp. $n\geq 15$, $n\geq 16$), the K3 surface contains a set of 13 (resp. 14,15) disjoint rational curves and this is possible if and only if it contains also 2 (resp. 3, 4) independent divisible classes made up of these curves (by Proposition \ref{prop: the setes Rm}), so if $n\geq 14$  (resp. $n\geq 15$, $n\geq 16$), then $NS(Y_n)$ is an overlattice of index at least $2^2$ (resp. $2^3$, $2^4$) of $L_n^{(1)}$.

The unique cases which remain open are the cases where $c$ is divisible, which can appear if $(n,a)=(9,2)$, $(13,2)$, $(14,3)$, $(17,2)$. We already discussed the divisibility of $c$ in case $n=9$, obtaining the lattice $L_9^{(2)}$. We excluded the case $n=17$. Let us consider the case $n=13$. In this case, since $(c-\sum_{i=1}^{12}r_i)/2$ is contained in $NS(Y_n)$, if $c/2$ is contained in $NS(Y_n)$, then we have that $\sum_{i=1}^{12}r_i/2\in NS(Y_n)$, but the sum of 12 disjoint rational curves can not be 2-divisible on a K3 surface. So we exclude also this case.
It remains only the case $(n,a)=(14,3)$. Since we have 13 disjoint rational curves, we know that there are exactly two 2-divisible subset of 8 disjoint rational curves, so obtain an overlattice of index $2^23$ of $L_n^{(1)}$. This lattice is a 2-elementary lattice of rank 14 length 8 and its discriminant form take values in $\Z$. But there exist no K3 surfaces whose N\'eron--Severi group is
characterized by these invariants (see \cite[Table 1]{NikFactosGroup}).\endproof
\begin{definition}\label{defi: the list of lattices}
Let us denote by $\mathcal{L}$ the following list of lattices: $$\mathcal{L}:=\left\{\begin{array}{llll}L_n^{(r)}&\mbox{ with }&r=1,2,4,6,8,16&\mbox{ and }\\ &\mbox{if }& r=1,& \mbox{ then }6\leq n\leq 12;\\ &\mbox{if }& r=2,& \mbox{ then }9\leq n\leq 13;\\ &\mbox{if }& r=4,& \mbox{ then }13\leq n\leq 14;\\ &\mbox{if }& r=8,& \mbox{ then }n=15;\\ &\mbox{if }& r=16,& \mbox{ then }n=16. \end{array}\right\}.$$
\end{definition}
So the Proposition \ref{prop: possible NS} can be rephrased as: if $Y_n$ is a K3 surface as in Condition \ref{condition}, then $NS(Y_n)\in\mathcal{L}$.

\subsection{The existence and the geometry of $Y_n$ and the existence of $X_n$}\label{subsec: existence h=-1}

In this section we prove the existence of the K3 surfaces $Y_n$ admitting the even set of curves required by Condition \ref{condition}. This immediately implies the existence of the surfaces $X_n$, double covers of $Y_n$ branched along $C\bigcup\cup_i R_i$.
In the following, for each value of $(n,r)$ such that $L_n^{(r)}\in \mathcal{L}$, we prove the existence of the surface $Y_n^{(r)}$ such that $NS(Y_n^{(r)})\simeq L_n^{(r)}$, we describe the geometry of the surface $Y_n^{(r)}$ and we prove that the classes $c$ and $r_i$ considered in Section \ref{subsec: possible NS} are supported of smooth irreducible curves. This allows to conclude the existence of the smooth double cover $X_n\ra Y_n^{(r)}$ as in Condition \ref{condition}. The results of this section (Propositions \ref{prop: Y_n with NS=Ln1}, \ref{prop: Y_n with NS=Ln2} and \ref{prop: Y_n with NS=Lni bog index}) together with Proposition \ref{prop: possible NS} imply the  following theorem, which concludes the classification of the K3 surfaces admitting a double cover as in Condition \ref{condition}:
\begin{theorem}\label{theorem: classification gC>2, h=-1 }
There exist K3 surfaces as in Condition \ref{condition} for every integer $n$ such that $6\leq n\leq 16$.
More precisely, a K3 surface $Y_n$ as in Condition \ref{condition} exists if and only if $NS(Y_n)\in\mathcal{L}$.
\end{theorem}

\subsubsection{The lattice $L_n^{(1)}$}\label{subsec: Ln1}

\begin{proposition}\label{prop: Y_n with NS=Ln1} If $6\leq n\leq 12$, then there exists a K3 surface $Y_n^{(1)}$ such that $NS(Y_n^{(1)})\simeq L_n^{(1)}$. 
Moreover one can assume that $c$ is the class of a pseudoample divisor on $Y_n^{(1)}$. 

If $n>6$, $\varphi_{|c|}(Y_n^{(1)})\subset \mathbb{P}^{n-4}$ has $n-1$ nodes, all these nodes are contained in the same hyperplane $H\simeq \mathbb{P}^{n-5}\subset\mathbb{P}^{n-4}$, they are all contained in a rational normal curve in $\mathbb{P}^{n-5}$ and the hyperplane section $H\cap\varphi_{|c|}(Y_n^{(1)})$ consists of this rational normal curve with multiplicity 2.

If $n=6$, the map $\varphi_{|c|}:Y_{6}^{(1)}\ra \mathbb{P}^2$ is a 2:1 cover branched along the union of a line and a smooth irreducible quintic.

The set of curves given by (the pull back of) a smooth hyperplane section of $\varphi_{|c|}(Y_n^{(1)})$ and the rational curves contracted by $\varphi_{|c|}$ is 2-divisible. 
\end{proposition}
\proof

{\bf Existence of $Y_n^{(1)}$.}  The rank of $L_n^{(1)}$ is $n$ and its length is $n-2$. By \cite[Theorem 1.14.4]{Nikbilinear}, a lattice $L$ admits a primitive embedding in $\Lambda_{K3}$ if $l(L)\leq 20-\rk(L)$, we conclude $L_n^{(1)}$ admits a primitive embedding in $\Lambda_{K3}$ if $n\leq 11$. 
In order to extend this result to the case $n=12$ we consider the discriminant group of $L_{12}^{(1)}$. It is $(\Z/14\Z)\times (\Z/2\Z)^9$ and it is generated by the vectors $d_1:=c/14+r_1/2$, $d_2:=(r_1+r_2)/2$, $d_3:=(r_2+r_3)/2$, $d_4:=(r_3+r_4+r_7+r_8)/2$, $d_5:=(r_4+r_5)/2$, $d_6:=(r_5+r_6)/2$, $d_7:=(r_8+r_9)/2$, $d_8:=(r_9+r_{10})/2$, $d_9:=(r_7+r_8+r_9+r_{10})/2$, $d_{10}:=(r_4+r_5+r_6+r_7)/2$. Since the sublattice generated by $d_9$ and $d_{10}$ is orthogonal to the one generated by $d_i$, $i=1,\ldots 8$ with respect to the discriminant form and the discriminant form computed on $d_9$ and $d_{10}$ is the same as the discriminant form of $U(2)$, \cite[Theorem 1.14.4]{Nikbilinear} guarantees that $L_{12}^{(1)}$ can be primitively embedded in $\Lambda_{K3}$.

By the Torelli theorem, there exists a K3 surface (and in fact a $(20-n)$-dimensional space of K3 surfaces) $Y_n^{(1)}$ such that $NS(Y_n^{(1)})\simeq L_n^{(1)}$. 

{\bf The divisor $c$ is pseudoample and without fixed components}. The class $c$ can be chosen to be a pseudoample class, since it has a positive square and we can assume that the chamber of the positive cone which corresponds to the ample cone has $c$ in its closure. 
By \cite[Theorem page 79]{R}, since $c$ is pseudoample, it has a fixed part if and only if it is $c=aE+\Gamma$ where $a\in\N$, $|E|$ is a free pencil and $\Gamma$ an irreducible $(-2)$-curve such that
$E\Gamma = 1$. Here we exclude this case if $n>6$. Indeed, if $c\simeq aE+\Gamma$, then $cE=(aE+\Gamma)E=1$. So there should exist a class $E$ such that $cE=1$. Each class in $NS(Y_n^{(1)})\simeq L_n^{(1)}$ is of the form $\alpha c+\sum_{j=1}^{n-1}\beta_j r_j$, where either $\alpha,\beta_j\in \Z$ or all the numbers $\alpha$, $\beta_j$, $j=1,\ldots n-1$ are non integer and contained in $\frac{1}{2}\Z$. So we are looking for $\alpha$, $\beta_j$, $j=1,\ldots, n-1$ such that $(\alpha c+\sum_{j=1}^{n-1}\beta_j r_j)c=(2n-10)\alpha=1$. This is impossible if $n>6$. Hence if $n>6$, $|c|$ has no a fixed part.

{\bf The linear system $|c|$ defines a map to $\mathbb{P}^{n-4}$, which is 1:1 provided that $n>6$.} Indeed, by \cite[Section (4.1)]{SD}, $\varphi_{|c|}$ is either 1:1 or 2:1 and by \cite[Theorem 5.2]{SD} it is 2:1 if and only if one of the following conditions hold: 
\begin{itemize}\item there exists a smooth curve $B$ of genus 2 such that $c\simeq 2B$ \item there exists a smooth curve $E$ of genus 1 such that $cE=2$.\end{itemize}
The existence of the curve $B$ would implies that $c/2\in NS(Y_n^{(1)})\simeq L_n^{(1)}$, which is easily proved to be false, since we have an explicit $\Z$ basis of $L_n^{(1)}$. The existence of the curve $E$ can be excluded in a similar way. Indeed if $E$ exists, its class is $\alpha c+\sum_{j=1}^{n-1}\beta_j r_j$, where either $\alpha,\beta_j\in \Z$ or all the numbers $\alpha$, $\beta_j$, $j=1,\ldots n-1$ are non integer and contained in $\frac{1}{2}\Z$.
Since $E$ should be a smooth curve of genus 1, $(\alpha c+\sum_{j=1}^{n-1}\beta_j r_j)^2=(2n-10)\alpha^2-2\sum_{j=1}^{n-1} \beta_j^2=0$ and since $cE=2$, $c(\alpha c+\sum_{j=1}^{n-1}\beta_j r_j)=2(n-5)\alpha=2$.
Hence we have two possibilities: either $n-5=\alpha=1$, but we are assuming $n>6$, or $n-5=2$ and $\alpha=\frac{1}{2}$. In the latter case $(\frac{1}{2} c+\sum_{j=1}^{6}\beta_j r_j)^2=1-2\sum_j \beta_j^2\leq 1-2\frac{6}{4}=-2$, so the class of $E$ does not satisfy $E^2=0$ and hence such a curve does not exist on $Y_n^{(1)}$.

{\bf The surface $\varphi_{|c|}(Y_n^{(1)})$ has $n-1$ nodes, all contained in the same hyperplane.} The orthogonal to $c$ in $L_n^{(1)}$, $c^{\perp_{L_n^{(1)}}}$ is spanned by the classes $r_i$. So the classes $r_i$ form a bases of the root lattice orthogonal to $c$. By \cite[Lemma 3.1]{B}, there exists a basis of the root lattice of $c^{\perp_{L_n^{(1)}}}$ which is supported on smooth irreducible rational curves. Thus we can assume that $r_i$ are the classes of smooth disjoint rational curves. Since $cr_i=0$, the curves $r_i$ are contracted to simple nodes on $\varphi_{|c|}(Y_n^{(1)})$.

We recall that the vector $L_n^{(1)}$ contains also the class $d_n:=(c-\sum_{j=1}^{n-1}r_j)/2$.
Since $d_n^2=-2$ and $d_nc=n-5>0$, the class $d_n$ is effective, and thus also $2d_n$ is effective.
But if the nodes $\varphi_{|c|}(Y_n^{(1)})\subset\mathbb{P}^{n-4}$ would not be contained in a hyperplane, the divisor $2d_n=c-r_1-\ldots -r_{n-1}$ is not effective, since it corresponds to the hyperplane passing through the nodes $\varphi_{|c|}(r_i)$. Thus all the nodes are on the same hyperplane, denoted by $H$. By construction the class of the hyperplane section $H\cap\varphi_{|c|}(Y_n^{(1)})$ is $2d_n$, where $d_n$ is an effective class. If the curve represented by $d_n$ would be reducible, then there exist two vectors $w_n, u_n\in L_n^{(1)}$, possibly coinciding, such that $d_n=w_n+u_n$, $u_n$ and $w_n$ represent effective divisors and do not contain the curves represented by $r_i$ as components. So we can assume that $w_n=\alpha c+\sum_{j=2}^{n}\beta_j r_j$, $u_n=\gamma c+\sum_{j=2}^{n}\delta_j r_j$ with $\alpha$, $\beta_j$, $\gamma$, $\delta_j\in \frac{1}{2}\Z$. By  $w_nc\geq 0$, $u_nc\geq 0$, we deduce that $\alpha\geq 0$, $\gamma\geq 0$; from $d_n=w_n+u_n$  $\alpha+\gamma=\frac{1}{2}$ we deduce that $\alpha=0$ and $\gamma=\frac{1}{2}$ (or viceversa). Similarly, by $w_nr_i\geq 0$, $u_nr_i\geq 0$, $d_n=w_n+u_n$ and $\gamma=\frac{1}{2}$, one deduces that $w_n=0$ and $d_n=u_n$. 
Hence $\varphi_{|c|}(d_n)\subset H\cap\varphi_{|c|}(Y_n^{(1)})$ is the class of a rational curve in $H\simeq \mathbb{P}^{n-5}$, of degree $n-5=cd_n$, i.e. it is a rational normal curve in $\mathbb{P}^{n-5}$. Since $d_nr_i=1$, the curve represented by $d_n$ passes through $n-1$ given points (and this determines such a curve in $\mathbb{P}^{n-5}$ uniquely) and the hyperplane section $H\cap\varphi_{|c|}(Y_n^{(1)})$ is this curve with multiplicity 2.

{\bf The case $n=6$} In this case $c^2=2$. Since there exists no class $E$ such that $cE=1$ and $E^2=0$, the class $c$ can not be expressed as $a E+\Gamma$ where $|E|$ is a free pencil, $\Gamma$ is a rational curve such that $E \Gamma=1$. Thus, by \cite[Theorem page 79]{R}, $|c|$ has no fixed part. Hence $\varphi_{|c|}:Y_n^{(1)}\ra \mathbb{P}^2$ is a $2:1$ cover, branched on a (possibly singular and reducible) sextic. The classes $r_j$ correspond to irreducible rational curves (as in case $n>6$) and they are contracted by $\varphi_{|c|}$. So the branch locus of the double cover $\varphi_{|c|}:Y_n^{(1)}\ra \mathbb{P}^2$ has 5 singular points, which are simple nodes. The class $d_2$ has self intersection $-2$ and $cd_2=1$, so $d_2$ is effective and $\varphi_{|c|}(d_2)$ is a line passing through the nodes. The pullback of the branch locus is represented by the class $3c-r_1-r_2-r_3-r_4-r_5$, which split into the sum of $d_2$ and $d_2+2c$. So the branch locus splits in the union of 2 curves, one of them is the line $\varphi_{|c|}(d_2)$ and the other is the quintic $\varphi_{|c|}(d_2+2c)$. The latter is irreducible because otherwise there should be other divisible classes in $L_6^{(1)}$ which represent the components of this quintic (alternatively one can observe that $NS(Y_{6}^{(1)})\simeq L_6^{(1)}$ is a 2-elementary lattice of rank 6 and length 4, so $Y_{6}^{(1)}$ admits a non symplectic involution acting trivially on $NS(Y_{6}^{(1)})$ which is exactly the cover involution of the map $\varphi_{|c|}:Y_{6}^{(1)}\ra\mathbb{P}^2$. The rank and the length of $NS(Y_{6}^{(1)})$ allows to compute the topology of the fixed locus of this involution, and those of the branch of the cover, which consists of the union of a smooth rational curve and a smooth curve of genus 6, cf. \cite{NikFactosGroup}). 

{\bf The K3 surface $Y_n^{(1)}$ admits a 2-divisible set and thus has a smooth cover $X_n$.} This is clear by what we said above: If $n>6$, let us consider a smooth hyperplane section of $\varphi_{|c|}(Y_n^{(1)})$ with a hyperplane which does not contain any nodes of $\varphi_{|c|}(Y_n^{(1)})$. It is a smooth irreducible curve on the K3 surface, denoted by $C$ and the class of $C$ is exactly the vector $c$. The curve $C$ is orthogonal to the smooth irreducible curves contracted by $\varphi_{|c|}$, whose classes are $r_j$. Since $d_n=(c-\sum_{i=1}^{n-1} r_i)/2\in NS(Y_n^{(1)})$, we have that $c+\sum_{i=1}^{n-1}r_i=2(d_n+\sum_{i=1}^{n-1}r_i)$ and thus  $c+\sum_{i=1}^{n-1}r_i$ is 2-divisible. This allows us to construct the double cover $X_n$ as in Condition \ref{condition}. The case $n=6$ is analogous except for the fact that the curve $C$ is a $2:1$ cover of the hyperplane section of projective space $|c|\simeq\mathbb{P}^2$.\endproof

\begin{proposition}\label{prop: projection of Y_n on P2}
If $NS(Y_n^{(1)})\simeq L_n^{(1)}$ and $n>6$, the projection of $\varphi_{|c|}(Y_n^{(1)})\subset \mathbb{P}^{n-4}$ from the $(n-7)$-dimensional projective space generated by $n-6$ of the nodes $\varphi_{|c|}(r_i)\subset Y_n^{(1)}$ exhibits $Y_n^{(1)}$ as double cover of $\mathbb{P}^2$ branched along a sextic which splits in a line and a smooth quintic. There are $n-6$ conics in $\mathbb{P}^2$ such that each of these conics intersects the branch reducible sextic with even multiplicity in each of their intersection points.
\end{proposition}
\proof Let us consider the map $\varphi_{|c|}:Y_n^{(1)}\ra\mathbb{P}^{n-4}$ and let us choose $n-6$ nodes in $\varphi_{|c|}(Y_n^{(1)})$. They span a linear space, $P_{n-7}$, of dimension $n-7$ and the image of the projection of $\mathbb{P}^{n-4}$ by $P_{n-7}$ is $\mathbb{P}^2$. Since the degree of $\varphi_{|c|}(Y_n^{(1)})$ is $2n-10$, the intersection of $\varphi_{|c|}(Y_n^{(1)})$ with a generic linear space of dimension $n-6$ which contains $P_{n-7}$ consists of $2n-10$ points, counted with their multiplicity. Since $P_{n-7}$ contains $n-6$ nodes of $\varphi_{|c|}(Y_n^{(1)})$, there are other $2=(2n-10)-2(n-6)$ points of intersection between $\varphi_{|c|}(Y_n^{(1)})$ and a generic linear space of dimension $n-6$ which contains $P_{n-7}$. So the projection of  $\varphi_{|c|}(Y_n^{(1)})$ by $P_{n-7}$ is a $2:1$ map $\varphi_{|c|}(Y_n^{(1)})\ra\mathbb{P}^2$. The image, under the projection from $P_{n-7}$, of the rational normal curve $H\cap \varphi_{|c|}(Y_n^{(1)})$ is in the branch locus of the map $2:1$ $\varphi_{|c|(Y_n^{(1)}\ra)}\ra\mathbb{P}^2$, since it has multiplicity 2. Without loss of generality we can assume that the $n-6$ nodes spanning $P_{n-7}$ are $\varphi_{|c|}(R_i)$, $i=1,\ldots, n-6$. Then the $2:1$ map $\varphi_{|c|}(Y_n^{(1)})\ra\mathbb{P}^2$ is induces on $Y_n^{(1)}$ by the linear system $|c-\sum_{i=1}^{n-6}r_i|$. We will denote by $h$ the divisor $c-\sum_{i=1}^{n-6}r_i$. The curves $R_{i}$, $n-5\leq i\leq n-1$ are contracted by the map $\varphi_{|h|}:Y_n^{(1)}\ra \mathbb{P}^2$ (indeed they correspond to the nodes of $\varphi_{|c|}(Y_n^{(1)})$ which are not contained in $P_{n-7}$). The class of the branch locus is $3h-\sum_{i=n-5}^{n-1} r_i=3(c-\sum_{i=1}^{n-6}r_i)-\sum_{i=n-5}^{n-1} r_i$. The branch curve is a reducible sextic with 5 nodes, the images of $R_i$ for $n-5\leq i\leq n-1$. It is reducible since the image of the rational normal curve $H\cap \varphi_{|c|}(Y_n^{(1)})$ is sent, by the projection to $\mathbb{P}^2$, to a component of the branch locus. Since its class is $d_n:=(c-\sum_{i=1}^{n-1}r_i)/2$ and $hd_n=(c-\sum_{i=2}^{n-5}r_i)d_n=(2n-10-2(n-6))/2=1$, one deduces that $\varphi_{|h|}(d_n)$ is a line in the branch locus. Thus the curve  $3h-\sum_{i=n-5}^{n-1}r_i-(c-\sum_{i=1}^{n-1}r_i)/2=3(c-\sum_{i=1}^{n-6}r_i)-\sum_{i=n-5}^{n-1}r_i-(c-\sum_{i=1}^{n-1}r_i)/2=(5h-\sum_{j=n-5}^{n-1} r_j)/2$ is a (possibly reducible) quintic in the branch locus. The points of intersection between the line $\varphi_{|h|}(d_n)$ and the quintic $\varphi_{|h|}((5h-\sum_{j=n-5}^{n-1} r_j)/2)$ are the 5 points $\varphi_{|h|}(r_i)$ with $i\geq n-5$. Since $((5h-\sum_{j=n-5}^{n-1} r_j)/2)^2=10$, the curve represented by this class has genus 6, which is indeed the genus of a smooth plane quintic. We now consider the image of the curves $R_i$ if $i\leq n-6$. In this case $hR_i=2$, thus $\varphi_{|h|}(R_i)$ is either a line or a plane conic, which splits in the double cover $\varphi_{|h|}:Y_n^{(1)}\ra\mathbb{P}^2$. The first case is impossible, since the pullback of a line not contained in the branch locus on $Y_n^{(1)}$ is linearly equivalent either to the class $h$ or to a class $h-R_j$, for $j\geq n-5$ but $R_i$ is linearly equivalent neither to $h$ nor to $h-R_j$. Thus $\varphi_{|h|}(R_i)$ is a plane conic, which splits in the double cover. Indeed both the classes $r_i$ and $2h-r_i$ are mapped on the same conic by $\varphi_{|h|}$. To recap: the map $\varphi_{|h|}:Y_n^{(1)}\ra\mathbb{P}^2$ is a 2:1 map whose branch locus has 5 singularities (in the points $\varphi_{|h|}(R_j)$, $n-5<j\leq n$) and is the union of a line (which is  $\varphi_{|h|}(d_n)$) and a smooth quintic ($\varphi_{|h|}((5h-\sum_{j=n-5}^{n-1} r_j)/2)$). There are $n-6$ conics which split in the double cover (which are $\varphi_{|h|}(R_i)$, for $i\leq n-6$). Since $R_iR_j=0$ if $i\neq j$, these conics never pass through the nodes of the branch locus.  Since they split in the double cover, their intersection with the branch locus has an even multiplicity in every intersection point and since they never pass through the singular points of the branch locus, this implies that they are tangent to the branch locus in every intersection points. Moreover, if $n>7$, two of these conics never intersect in points which are on the branch locus (otherwise their inverse image on $Y_n^{(1)}$ should intersect, but this never happens since $R_iR_j=0$ if $i\neq j$).\endproof

\begin{rem}\label{rem: 6 conics tg to a quinti and a line}{\rm Let us consider the K3 surface $Y_{12}^{(1)}$ such that $NS(Y_{12}^{(1)})\simeq L_{12}^{(1)}$. Then the map $\varphi_{|h|}$ defined in Proposition \ref{prop: projection of Y_n on P2} exhibits $Y_{12}^{(1)}$ as double cover of $\mathbb{P}^2$ branched along the union of a quintic $Q$ and a line $L$. Moreover there are 6 conics $F_i$, $i=1,\ldots 6$, with special properties with respect to $Q$ and $L$, indeed  they are such that:\begin{itemize}\item $F_i\cap Q\cap L=\emptyset$ for every $i$; \item the multiplicity of intersection of $F_i$ and $Q$ is even in each of their intersection points (in particular $F_i$ is tangent to $Q$ in every intersection points); \item the multiplicity of intersection of $F_i$ and $L$ is even in each of their intersection points (in particular $F_i$ is tangent to $L$ in every intersection points);\item an intersection point of the  two conics $F_i$ and $F_j$ is never contained in $Q\cup L$.\end{itemize}
Each of these conics is an odd theta characteristic of $Q$, but they also have the peculiarity that all of them are tangent to the same line.

We recall  that it is always possible to find tangent conics to a smooth quintic, since as we noticed, these correspond to the odd theta characteristics of the quintic. But in order to find the previous configurations of $Q$, $L$ and $F_i$ one has to fix a quintic which is such that there exists 6 odd theta characteristics which are associated to conics with a commune tangent line. The existence of this configuration can be proved exactly by the existence of the K3 surface $Y_{12}^{(1)}$. Moreover, since $Y_{12}^{(1)}$ lives in an 8 dimensional space, we expect that 8 is also the dimension of the space of quintics with this property. In Corollary \ref{cor: existence theta characteristics} we improve the bound on the number of conics with the same properties as the $F_i$'s.}\end{rem}

\begin{proposition}
Let $W_n$ be a K3 surface which admits a model $\varphi(W_n)\subset\mathbb{P}^{n-4}$ whose singularities are at worst simple double points. If there is a hyperplane $\Pi\subset\mathbb{P}^{n-4}$ such that $\varphi(W_n)\cap\Pi$ is a rational normal curve in $\Pi\simeq \mathbb{P}^{n-5}$ with multiplicity 2, then $\varphi(W_n)$ has $n-1$ nodes, all contained in the same hyperplane $\Pi$ of $\mathbb{P}^{n-4}$.

If the Picard number of $W_n$ is $n$, then $W_n$ admits a 2:1 cover as in Condition \ref{condition}. 
\end{proposition}
\proof Let us denote by $H$ the pseudoample divisor such that $\varphi_{|H|}:W_n\ra\mathbb{P}^{n-4}$. Then $H^2=2n-10$. Each hyperplane section of $\varphi_{|H|}$ is represented either by the class $H$ or by the class $H-\sum_{i=1}^{f} M_i$ where $M_i$ are classes in $NS(W_n)$ which represent smooth rational curves contracted by $H$. In this case we are considering the hyperplane section of $\varphi_{|H|}(W_n)$ passing through $f$ nodes, which are $\varphi_{|H|}(M_i)$. Since the curves $M_i$ are contracted by $\varphi_{|H|}$, $HM_i=0$ and since the singularities of $\varphi_{|H|}(W_n)$ are at worst double points, $M_iM_j=0$ if $i\neq j$. Moreover $M_i^2=-2$, because $M_i$ are smooth rational curves. Let us assume that a hyperplane section of $\varphi_{|H|}(W_n)$ is a curve $D$ counted with multiplicity 2. Thus the class of $2D$ is a hyperplane section, so either $2D\simeq H$ or $2D\simeq H-\sum_{i=1}^{f} M_i$. By hypothesis $D$ is a smooth rational curves on $W_n$, so $D^2=-2$, so either $(H/2)^2=-2$ or  $\left(\left(H-\sum_{i=1}^{f}M_i\right)\right)^2=-2$. The condition $(H/2)^2=-2$ implies $H^2=-8$, which is not possible. The condition $\left(\left(H-\sum_{i=1}^{f}M_i\right)/2\right)^2=-2$ implies $n-5-f=-4$, i.e. $f=n-1$. So there exist $n-1$ nodes on $\varphi_{|H|}(W_n)$ and a hyperplane passing through all these nodes which cuts a rational curve on $\varphi_{|H|}(W_n)$ with multiplicity 2, whose class is $D$. The degree of the curve represented by $D$ is $DH=n-5$, hence $D$ is a rational normal curve in $\mathbb{P}^{n-5}$.
By the assumptions the classes $H$ and $M_i$ form a 2-divisible set of disjoint curves, so there exists a smooth double cover of $W_n$ and it is as in Condition \ref{condition}, since the genus of a generic curve in $|H|$ is $H^2/2+1=n-4$. 
\endproof

\subsubsection{The lattice $L_n^{(2)}$}\label{subsec: Ln2}

\begin{proposition}\label{prop: Y_n with NS=Ln2} If $9\leq n\leq 13$, then there exists a K3 surface $Y^{(2)}_n$ such that $NS(Y_n^{(2)})\simeq L_n^{(2)}$. 
Moreover one can assume that $c$ is the class of a pseudoample divisor on $Y_n^{(2)}$ and if $n>9$, $\varphi_{|c|}(Y_n^{(2)})\subset \mathbb{P}^{n-4}$ has $n-1$ nodes. 

All these nodes are contained in the same hyperplane $H\simeq \mathbb{P}^{n-5}\subset\mathbb{P}^{n-4}$, they are all contained in a rational normal curve in $\mathbb{P}^{n-5}$ and the hyperplane section $H\cap\varphi_{|c|}(Y_n^{(2)})$ consists of this rational normal curve with multiplicity 2. 

If $n=9$, the map $\varphi_{|c/2|}:Y_9^{(2)}\ra \mathbb{P}^2$ is a 2:1 cover branched along the union of a conic and a quartic. Denoted by $V_2:\mathbb{P}^2\ra\mathbb{P}^5$ the Veronese map, $\varphi_{|c|}:Y_9^{(2)}\ra \mathbb{P}^5$ is the composition $ V_2\circ \varphi_{|c/2|}$.

In all the cases the K3 surface $Y_n^{(2)}$ admits a smooth double cover $X_n^{(2)}\ra Y_n^{(2)}$ as in Condition \ref{condition}. 
\end{proposition}
\proof
The proof is very similar to the one of Proposition \ref{prop: Y_n with NS=Ln1}. In particular the proof of  the existence of the K3 surfaces $Y_n^{(2)}\simeq L_n^{(2)}$ is analogous and we omit it.

As in proof of Proposition \ref{prop: Y_n with NS=Ln1}, one can assume that $c$ is pseudoample and using \cite[Theorem page 79]{R} and \cite[Section (4.1) and Theorem 5.2]{SD} one proves that $|c|$ has no fixed part and defines a $1:1$ map if $n>9$. In particular we observe that $\varphi_{|c|}$ is $2:1$ if there exists a smooth curve $B$ of genus 2 such that $c\simeq 2B$.
The existence of the curve $B$ would imply that $c/2\in NS(Y_n^{(2)})\simeq L_n^{(2)}$, which is easily proved to be false if $n>9$ (we observe that on the other hand this is true if $n=9$). 

The root lattice of  $c^{\perp_{L_n^{(2)}}}$ is spanned by the classes $r_i$.  By \cite[Lemma 3.1]{B}, there exists a basis of the root lattice $c^{\perp_{L_n^{(2)}}}$ which is supported on smooth irreducible rational curves. Thus we can assume that $r_i$ are the classes of smooth disjoint rational curves. Since $cr_i=0$, the curves $r_i$ are contracted to simple nodes on $\varphi_{|c|}(Y_n^{(2)})$. As in proof of Proposition \ref{prop: Y_n with NS=Ln1} one has that the class $d_n:=(c-\sum_{j=2}^nr_j)/2\in L_n^{(2)}$ is effective and so all the nodes are on the same hyperplane, denoted by $H$. By construction the class of the hyperplane section $H\cap\varphi_{|c|}(Y_n^{(2)})$ is $2d_n$, where $d_n$ is an effective class. If the curve represented by $d_n$ would be reducible, then there exists two classes $u$ and $v$ such that $u+v\simeq d_n$, $uc>0$, $vc>0$. Since all the elements in $L_n^{(2)}$ are linear combinations of $c,r_i$, $i=1,\ldots, n-1$ with coefficients in $\frac{1}{2}\Z$, $u:=\alpha c+\sum_{i=1}^{n-1}\beta r_i$ and $v:=a c+\sum_{i=1}^{n-1}b r_i$ with $\alpha,\beta,a,b\in \frac{1}{2}\Z$ . By $uc>0$, it follows that $\alpha>0$ and by $vc>0$, it follows that $a>0$. By $u+v=d_n$, it follows that $\alpha+a=\frac{1}{2}$, which is impossible, since both $\alpha$ and $a$ are contained in $\frac{1}{2}\Z$ and are positive.

{\bf The case $n=9$.} In this case $c^2=8$. Since the class $\sum_{i=1}^8r_i/2$ is contained in $NS(Y_n^{(2)})$, also the class $d_n+\sum_{i=1}^8r_i/2=c/2$ is contained in $NS(Y_n^{(2)})$. There exists no a class $E$ such that $(c/2)E=1$ and $E^2=0$, so the class $c/2$ can not be expressed as $a E+\Gamma$ where $|E|$ is a free pencil, $\Gamma$ is a rational curve such that $E \Gamma=1$. Thus, by \cite[Theorem page 79]{R}, $|c/2|$ has no a fixed part. Hence $\varphi_{|c/2|}:Y_n\ra \mathbb{P}^2$ is a $2:1$ cover, branched on a (possibly singular and reducible) sextic. The classes $r_j$ correspond to irreducible rational curves (as in case $n>9$) and they are contracted by $\varphi_{|c/2|}$. So the branch locus of the double cover $\varphi_{|c/2|}:Y_9^{(2)}\ra \mathbb{P}^2$ has 8 singular points, which are simple nodes. The class $d_9$ has self intersection $-2$ and $(c/2)d_9=2$, so $d_9$ is effective and $\varphi_{|c/2|}(d_9)$ is a conic passing through the points $\varphi_{|c/2|}(r_i)$. The pulllback of the branch locus is represented by the class $3(c/2)-r_1-r_2-r_3-r_4-r_5-r_6-r_7-r_8$, which splits into the sum of $d_9$ and $d_9+c/2$. So the branch locus splits in the union of 2 curves, one of them is the conic $\varphi_{|c/2|}(d_9)$ and the other is $\varphi_{|c/2|}(d_9+c/2)$. Since $(c/2)(d_9+c/2)=4$ and $(d_9+c/2)^2=4$, the curve  $\varphi_{|c/2|}(d_9+c/2)$ is an irreducible smooth plane quartic curve. As in proof of Proposition \ref{prop: Y_n with NS=Ln1}, one can also observes that $Y_9^{(2)}$ admits a non symplectic involution whose fixed locus is a rational curve and a curve of genus 3.
The divisor $c$ is hyperelliptic, by \cite[Theorem 5.2]{SD}, since there exists a curve $B$ of genus 2, whose class is $c/2$, such that $c\simeq\mathcal{O}(B)$. The map induced by $|c|$ is described in \cite[Proposition 5.6]{SD}: $\varphi_{|c|}:Y_n^{(2)}\ra\mathbb{P}^5$, is the composition of the map $\varphi_{|c/2|}$ and the Veronese embedding $\mathbb{P}^2\ra\mathbb{P}^5$.
 
{\bf There exists the smooth double cover $X_n^{(2)}$.} This is clear by what we said above: If $n>9$, let us consider a smooth hyperplane section of $\varphi_{|c|}(Y_n^{(2)})$ with a hyperplane which does not contains any nodes of $\varphi_{|c|}(Y_n^{(2)})$, $n>9$. It is a smooth irreducible curve on the K3 surface, denoted by $C$ and the class of $C$ is exactly the vector $c$. The curve $C$ is orthogonal to the smooth irreducible curves contracted by $\varphi_{|c|}$, whose classes are $r_j$. Since $d_n=(c-\sum_{i=1}^{n-1} r_i)/2\in NS(Y_n^{(2)})$, we have that $c+\sum_{i=1}^{n-1}r_i=2(d_n+\sum_{i=1}^{n-1}r_i)$ and thus $c+\sum_{i=1}^{n-1}r_i$ is 2-divisible. Similarly, if $n=9$, we denote by $C$ a curve in $|c|$, which is a smooth curve of genus 5, and by $R_i$ the curve contracted by $\varphi_{|c|}$ and we obtain that $\{C, R_i\}$ is a divisible set and by construction the double cover branched on $C\cup_i R_i$ satisfies the Condition \ref{condition}. \endproof

\begin{proposition}\label{prop: projection of Y_n(2) on P2}
If $NS(Y_n^{(2)})\simeq L_n^{(2)}$ and $n>9$, then there exists a choice of $n-6$ nodes $\varphi_{|c|}(R_i)\subset \varphi_{|c|}(Y_n^{(2)})$ such that  the projection of $\varphi_{|c|}(Y_n^{(2)})\subset \mathbb{P}^{n-4}$ from the $(n-7)$-dimensional projective space generated by these nodes exhibits $Y_n^{(2)}$ as double cover of $\mathbb{P}^2$ branched along a sextic which is the union of a line and a smooth quintic. There are $n-6$ conics such that each of these intersects the branch (reducible) sextic with even multiplicity in each intersection point.

For $n>9$, there also exists a choice of $n-9$ nodes $\varphi_{|c|}(R_i)\subset \varphi_{|c|}(Y_n^{(2)})$ such that the projection of $\varphi_{|c|}(Y_n^{(2)})\subset \mathbb{P}^{n-4}$ from the $(n-10)$-dimensional projective space generated by these nodes exhibits $Y_n^{(2)}$ as double cover of $V_2(\mathbb{P}^2)\subset \mathbb{P}^5$, where $V_2$ is the Veronese embedding $\mathbb{P}^2\ra\mathbb{P}^5$. Indeed the surface $Y_n^{(2)}$ admits a $2:1$ map to $\mathbb{P}^2$ branched along the union of a conic and a quartic. There are $n-9$ lines which intersect the branch locus with even multiplicity in the intersection points and which do not pass through the nodes of the branch locus.
\end{proposition}

\proof The proof is very similar to the one of Proposition \ref{prop: projection of Y_n on P2} and we omit it. We only make some useful remarks: one can assume that the nodes of $\varphi_{|c|}(Y_n^{(2)})$ which form a 2-divisible set are $\varphi_{|c|}(R_i)$ for $i=1,\ldots 8$. Then the projection from the linear space generated by the nodes $\varphi_{|c|}(R_i)$ for $i=1,\ldots n-6$ contracts exactly 5 rational curves of $Y_n^{(2)}$, i.e. the curves $R_i$ for $i=n-5,\ldots,n-1$. Indeed one can explicitly determine the number of $(-2)$-classes orthogonal to $c-\sum_{i=1}^{n-6} r_i$ and observe that it is 5, up to the choice of a sign. This implies the smoothness of the quintic which is a component of the branch locus of the $2:1$ map to $\mathbb{P}^2$ in the statement, indeed the branch locus has exactly 5 singular points which are the intersection points of the line (image of the class $d_n$) and the quintic.

In the second statement one has to observe that the divisor $d_n+(\sum_{i=1}^8r_i)/2=(c-\sum_{j=9}^{n-1} r_j)/2$ is contained in the N\'eron Severi group, hence the divisor $c-\sum_{j=9}^{n-1} r_j$ is 2-divisible. Since $c-\sum_{j=9}^{n-1} r_j$ corresponds to the projection of $\varphi_{|c|}(Y_n^{(2)})\subset \mathbb{P}^{n-4}$ from the $(n-10)$-dimensional projective space generated by $n-9$ of the nodes $\varphi_{|c|}(R_i)\subset \varphi_{|c|}(Y_n)$ with $n-9<i\leq n$ to $\mathbb{P}^2$, one obtains that the map $\varphi_{|c-\sum_{j=9}^{n-1} r_j|}$ is the composition of $\varphi_{|(c-\sum_{j=9}^{n-1} r_j)/2|}$ and the Veronese embedding $V_2$, by \cite[Proposition 5.6]{SD}. The irreducibility of the conic and the quartic in the branch locus of the double cover of $\mathbb{P}^2$ can be proved as before. The lines which are tangent to the branch locus are given by $\varphi_{|(c-\sum_{j=9}^{n-1} r_j)/2|}(R_k)$ for $k\geq 9$.
\endproof

\begin{rem}\label{rem: 4 lines tg to a quartic and a conic}{\rm Let us consider the K3 surface $Y_{13}^{(2)}$ such that $NS(Y_{13}^{(2)})\simeq L_{13}^{(2)}$. Then the map $\varphi_{\left|\left(c-\sum_{i=9}^{12}r_i\right)/2\right|}$ defined in Proposition \ref{prop: projection of Y_n(2) on P2} exhibits $Y_{13}^{(2)}$ as double cover of $\mathbb{P}^2$ branched along the union of a quartic $R$ and a conic $B$. Moreover there are 4 lines $G_i$, $i=1,\ldots, 4$, with special properties with respect to $R$ and $B$, indeed  they are such that:\begin{itemize}\item $G_i\cap R\cap B=\emptyset$ for every $i$; \item the multiplicity of intersection of $G_i$ and $R$ is even in each of their intersection points (in particular $G_i$ is tangent to $R$ in every intersection points); \item the multiplicity of intersection of $G_i$ and $B$ is even in each of their intersection points (in particular $G_i$ is tangent to $B$ in every intersection points);\item the intersection point among two lines $G_i$ and $G_j$ are never contained in $R\cup B$.\end{itemize}
Each of these lines is an odd theta characteristic of $R$. This result is similar to the one of Remark \ref{rem: 6 conics tg to a quinti and a line}, but in this case the existence of the configuration of curves $R,B$ and $G_i$ is a direct consequence of the existence of 4 bitangent lines for the quartic $R$. Indeed once one fixes $R$ and 4 bitangent lines $G_i$ to $R$, there is a 1-dimensional family of conics which are tangent to all these lines and thus it is surely possible to choose $B$ in this family (one possibly has to exclude the conics passing through the intersection points between $R$ and $G_i$).
In Corollary \ref{cor: existence theta characteristics} we improve the bound on the number of lines with the same properties of the $G_i$'s.}\end{rem}

\subsubsection{The lattices $L_n^{(r)}$, $r>2$}
In this section we state results on the lattices $L_n^{(r)}$ and on the K3 surfaces $Y_n^{(r)}$ such that $NS(Y_n^{(r)})\simeq L_n^{(r)}$ in case $r>2$. The results and the proofs are analogously to the ones described for the cases $r=1$ and $r=2$ in Sections \ref{subsec: Ln1} and \ref{subsec: Ln2} respectively. 

\begin{proposition}\label{prop: Y_n with NS=Lni bog index} Let $n=13,14$ (resp. $n=15$, $n=16$), then there exists a K3 surface $Y^{(4)}_n$ (resp. $Y_{15}^{(8)}$, $Y_{16}^{(16)}$) such that $NS(Y_n^{(r)})\simeq L_n^{(r)}$, $r=4$ (resp. $r=8$, $r=16$). 
Moreover, one can assume that $c$ is the class of a pseudoample divisor on $Y_n^{(r)}$ and $\varphi_{|c|}(Y_n^{(r)})\subset \mathbb{P}^{n-4}$ has $n-1$ nodes. 

All these nodes are contained in the same hyperplane $H\simeq \mathbb{P}^{n-5}\subset\mathbb{P}^{n-4}$, they are all contained in a rational normal curve in $\mathbb{P}^{n-5}$ and the hyperplane section $H\cap\varphi_{|c|}(Y_n^{(r)})$ consists of this rational normal curve with multiplicity 2. 

The K3 surface $Y_n^{(r)}$ admits a smooth double cover $X_n^{(r)}\ra Y_n^{(r)}$ as in Condition \ref{condition} for $r=4,8,16$. 
\end{proposition}

\begin{proposition}\label{prop: projection of Y_n(i) on P2 bog i}
If $NS(Y_n^{(r)})\simeq L_n^{(r)}$ $r=4,8,16$, then there exists a choice of $n-6$ nodes $\varphi_{|c|}(R_i)\subset \varphi_{|c}|(Y_n^{(r)})$  such that  the projection of $\varphi_{|c|}(Y_n^{(r)})\subset \mathbb{P}^{n-4}$ from the $(n-7)$-dimensional projective space spanned by these nodes exhibits $Y_n^{(r)}$ as double cover of $\mathbb{P}^2$ branched along a sextic which is the union of a line and a smooth quintic. There are $n-6$ conics such that each of these conics intersects the branch reducible sextic with even multiplicity in each intersection point.

There exists a choice of $n-9$ nodes $\varphi_{|c|}(R_i)\subset \varphi_{|c|}(Y_n^{(r)})$ such that the projection of $\varphi_{|c|}(Y_n^{(r)})\subset \mathbb{P}^{n-4}$ from the $(n-10)$-dimensional projective space generated by these nodes exhibits $Y_n^{(r)}$ as double cover of $V_2(\mathbb{P}^2)\subset \mathbb{P}^5$, where $V_2$ is the Veronese embedding $\mathbb{P}^2\ra\mathbb{P}^5$. Indeed the surface $Y_n^{(r)}$ admits a $2:1$ map to $\mathbb{P}^2$ branched along the union of a conic and a quartic. There are $n-9$ lines which intersect the branch locus with even multiplicity in the intersection points and which do not pass through the nodes of the branch locus.
\end{proposition}

\begin{corollary}\label{cor: existence theta characteristics}
There exist a quintic $Q$, a line $L$ and at least 9 different conics $F_i$, $i=1,\ldots, 9$ in $\mathbb{P}^2$ such that:\begin{itemize}\item $F_i\cap Q\cap L=\emptyset$ for every $i$; \item the multiplicity of intersection of $F_i$ and $Q$ is even in each of their intersection points (in particular $F_i$ is tangent to $Q$ in every intersection points); \item the multiplicity of intersection of $F_i$ and $L$ is even in each of their intersection points (in particular $F_i$ is tangent to $L$ in every intersection points);\item the intersection points among two conics $F_i$ and $F_j$ are never contained in $Q\cup L$.\end{itemize}
Each of these conics cuts an odd theta characteristic of $Q$. 

There exist a quartic $R$, a conic $B$ and at least 7 different lines $G_i$, $i=1,\ldots, 7$ in $\mathbb{P}^2$ such that:\begin{itemize}\item $G_i\cap R\cap B=\emptyset$ for every $i$; \item the multiplicity of intersection of $G_i$ and $R$ is even in each of their intersection points (in particular $G_i$ is tangent to $R$ in every intersection points); \item the multiplicity of intersection of $G_i$ and $B$ is even in each of their intersection points (in particular $G_i$ is tangent to $B$ in every intersection points);\item the intersection points among two conics $G_i$ and $G_j$ are never contained in $R\cup B$.\end{itemize}
Each of these lines cuts an odd theta characteristic of $R$. 
\end{corollary}
\proof The existence of these configurations follows directly by the existence of the K3 surface $Y_{16}^{(16)}$.\endproof

\begin{corollary}\label{cor: D is irreducible all cases}
With the same notation of Proposition \ref{prop: the effective class D}, the effective divisor $D$ on $Y_n$ such that $C\simeq 2D+\sum_{i=1}^{n-1}R_i$ represents a smooth irreducible rational curve if $NS(Y_n^{(r)})\simeq L_n^{(r)}$, for every admissible pair $(n,r)$.
\end{corollary}
\proof The class of the effective divisor $D$ coincides with the class $d_n:=c-\sum_{i=1}^n r_i$ which is sent by the map associated to a certain divisor, to a line $l$ in $\mathbb{P}^2$. In the (singular) double cover of $\mathbb{P}^2$ branched along the line $l$ and a smooth quintic, the pull back of the line $l$ is still irreducible. Then one blows up the 5 singular points of the double cover: the class $d_n$ is the class of the strict transform of the pull back of the line $l$, so it is the strict transform of an irreducible curve. We deduce that $D$ is irreducible on $Y_n^{(r)}$.\endproof

\subsection{The surfaces $X_{6}^{(1)}$ and $X_9^{(2)}$, bidouble covers of $\mathbb{P}^2$}\label{subsec: bidouble covers P2}

In this section we consider two interesting examples: the smooth double covers of the K3 surfaces $Y_{6}^{(1)}$ and $Y_9^{(2)}$. The peculiarity of these examples is that the N\'eron--Severi groups of both the surfaces $Y_{6}^{(1)}$ and $Y_9^{(2)}$ are 2-elementary, so each of these surfaces admits a non-symplectic involution which acts trivially on the N\'eron--Severi group. The main result is that this involution lifts to an involution of the double cover surfaces, $X_{6}^{(1)}$ and $X_9^{(2)}$. Indeed we will show that the surfaces $X_{6}^{(1)}$ and $X_9^{(2)}$ are bidouble covers of $\mathbb{P}^2$. This provides a very clear geometric description of the surfaces $X_{6}^{(1)}$ and $X_9^{(2)}$ and has an interesting consequence for the Hodge structure of their second cohmology groups.

\subsubsection{The surfaces $Y_{6}^{(1)}$ and $X_{6}^{(1)}$}\label{subsec: Y6 and X6}
Let $Y_{6}^{(1)}$ be a K3 surface such that $NS(Y_{6}^{(1)})\simeq L_6^{(1)}$. 
With the same notation as Proposition \ref{prop: Y_n with NS=Ln1}, we have that $\varphi_{|c|}:Y_{6}^{(1)}\ra\mathbb{P}^2$ is a double cover branched along the line $l:=\varphi_{|c|}(d_6)$ and the quintic $q:=\varphi_{|c|}\left((5c-\sum_{i=1}^5r_i)/2\right)$. We denote by $\iota_Y$ the cover involution on $Y_{6}^{(1)}$. 
By Proposition \ref{prop: Y_n with NS=Ln1}, $Y_{6}^{(1)}$ admits a double cover $X_{6}^{(1)}\ra Y_{6}^{(1)}$ which satisfies Condition \ref{condition}. 
Let us now consider the curve $C\subset Y_{6}^{(1)}$ which is the curve of genus 2 contained in the branch locus of the double cover $X_{6}^{(1)}\ra Y_{6}^{(1)}$. It is sent to a line by the map associated to its linear system, i.e. by $\varphi_{|c|}$. We denote by $m$ the line $\varphi_{|c|}(C)\subset \mathbb{P}^2$. We observe that $m\neq l$, and indeed $l$ is in the branch locus of the double cover $\varphi_{|c|}:Y_{6}^{(1)}\ra \mathbb{P}^2$ and $m$ is not.

\begin{proposition}
Let $S_6^{(1)}$ be the bidouble cover of $\mathbb{P}^2$ branched on $m\cup l\cup q$ with multiplicity 2 on each component. Then the surface $X_{6}^{(1)}$ is birational to $S_6^{(1)}$.
\end{proposition}
\proof The surface $S_6^{(1)}$, a bidouble cover of $\mathbb{P}^2$, is uniquely determined by the data $D_i$ $i=1,2,3$ and $L_i$, $i=1,2,3$ as in \cite[Proposition 2.3]{C}, where $D_i$ are the effective divisor corresponding to the branch of the three double covers of $\mathbb{P}^2$ and $2L_i\equiv D_k+D_j$. In our situation the curves in the branch locus are $l$, $m$ and $q$ and their classes in $\Pic(\mathbb{P}^2)\simeq \Z H$ are respectively $H$, $H$ and $5H$. Assuming $D_1\simeq [m]=H$, $D_2\simeq [l]=H$, $D_3\simeq [q]=5H$ we have $L_1\equiv 3H$, $L_2\equiv 3H$ and $L_3\equiv H$. These data determine $S_6^{(1)}$ uniquely. Moreover $S_6^{(1)}$ is a smooth surface, since $l$, $m$ and $q$ are smooth and have simple normal crossings.   
Let us denote by $Y_6^{(1)'}$ the double cover of $\mathbb{P}^2$ branched along $l\cup q$. The surface $Y_6^{(1)'}$ is a singular model of the K3 surface $Y_{6}^{(1)}$, obtained by the contraction of the 5 rational curves $R_i$. 
By construction $S_6^{(1)}$ is the double cover of $Y_6^{(1)'}$ branched along the 5 singular points of $Y_6^{(1)'}$ and the pull back of $m$ on $Y_6^{(1)'}$. Since the singular points of $Y_6^{(1)'}$ correspond to the 5 rational curves $R_i$ on $Y_{6}^{(1)}$ and the pull back of $m$ on $Y_6^{(1)'}$ is the image of $C\subset Y_{6}^{(1)}$ under the contraction map $Y_{6}^{(1)}\ra Y_6^{(1)'}$, $S_6^{(1)}$ is birational to the double cover of $Y_{6}^{(1)}$ branched along $C$ and the curves $R_i$, so it is birational to $X_{6}^{(1)}$.
\endproof

By construction, we have the following commutative diagram:
$$\xymatrix{&S_6^{(1)}\ar[dr]\ar[d]\ar[dl]&\\
Y_6^{(1)'}\ar[dr]&W_6^{(1)'}\ar[d]&V_6^{(1)'}\ar[dl]\\
&\mathbb{P}^2&}$$
where $W_6^{(1)'}$ is the double cover of $\mathbb{P}^2$ branched along $m\cup q$ and $V_6^{(1)'}$ is the double cover of $\mathbb{P}^2$ branched along $l\cup m$. The minimal model $W_6^{(1)}$ of $W_6^{(1)'}$ is a K3 surface (since $W_6^{(1)'}$ is a double cover of $\mathbb{P}^2$ branched on a nodal sextic) and each minimal model $V_6^{(1)}$ of $V_6^{(1)'}$ is a smooth rational surface. In particular $h^{2,0}(W_6^{(1)})=1$ and $h^{2,0}(V_6^{(1)})=0$.

The surface $S_6^{(1)}$ admits three commuting involutions, $\sigma_1$, $\sigma_2$ and $\sigma_3$, such that $S_6^{(1)}/\sigma_1\simeq Y_6^{(1)'}$, $S_6^{(1)}/\sigma_2\simeq W_6^{(1)'}$, $S_6^{(1)}/\sigma_3\simeq V_6^{(1)'}$. So the Hodge structure of $H^2(S_6^{(1)},\Q)$ splits in eigenspaces relative to the actions of $\sigma_1$ and $\sigma_2$. In particular we will denote by $H^2(S_6^{(1)},\Q)_{(\epsilon, \eta)}$, with $\{\epsilon, \eta\}\in \{+1,-1\}$ the intersection of the eigenspace of the value $\epsilon$ for $\sigma_1$ and the eigenspace of the value $\eta$ for $\sigma_2$.
So we have the following decomposition:
$$H^2(S_6^{(1)},\Q)=H^2(S_6^{(1)},\Q)_{(+1,+1)}\oplus H^2(S_6^{(1)},\Q)_{(-1,+1)}\oplus H^2(S_6^{(1)},\Q)_{(+1,-1)}\oplus H^2(S_6^{(1)},\Q)_{(-1,-1)}.$$

The Hodge structure on $H^2(S_6^{(1)},\Q)$ induces a Hodge structure on each eigenspace $H^2(S_6^{(1)},\Q)_{\epsilon, \eta}$. Since $S_6^{(1)}/\langle\sigma_1,\sigma_2\rangle$ is $\mathbb{P}^2$ and $h^{2,0}(\mathbb{P}^2)=0$, we have that $H^{2,0}(S_6^{(1)})_{(+1,+1)}$ is trivial. The same holds for $H^{2,0}(S_6^{(1)})_{(-1,-1)}$, which corresponds to the subspace invariant for $\sigma_3$, and $S_6^{(1)}/\sigma_3$ is a rational surface. The Hodge structure of both $H^2(S_6^{(1)},\Q)_{(-1,+1)}$ and $ H^2(S_6^{(1)},\Q)_{(+1,-1)}$ is of K3 type: one is induced by the Hodge structure of $H^2(Y_{6}^{(1)},\Z)$, and the other is induced by the Hodge structure of $H^2(W_6^{(1)},\Z)$.

Since the transcendental part of the second cohomology group of a surface is not influenced by birational transformations, we deduce that on $T_{X_{6}^{(1)}}$ we have two sub--Hodge structures of K3 type and each of them is induced by the Hodge structure of the transcendental part of the second cohomology group of a K3 surface. 

\begin{rem}{\rm Generically the N\'eron--Severi groups of the K3 surfaces $W_6^{(1)}$ and $Y_{6}^{(1)}$ are isometric, since both $Y_{6}^{(1)}$ and $W_6^{(1)}$ are the minimal resolution of the double cover of $\mathbb{P}^2$ branched along a line and a smooth quintic. However, the surfaces $W_6^{(1)}$ and $Y_6^{(1)}$ are not isomorphic, in general, since the specialization of one of them (assuming for example that there exists a plane curve which meets $m\cup q$ always with even multiplicity) does not imply the same specialization for the other. Thus, under certain choice, one can obtain different N\'eron--Severi groups for $W_6^{(1)}$ and $Y_{6}^{(1)}$.}\end{rem}

\begin{rem}{\rm
We underline that our construction is strongly related to the presence of the involution $\iota_Y$ on $Y_{6}^{(1)}$. Indeed the involution induced by $\sigma_2$ (and also by $\sigma_3$) on $Y_6^{(1)'}$ is the same involution induced by $\iota_Y\in$Aut$(Y_{6}^{(1)})$ on the contracted surface $Y_6^{(1)'}$. Equivalently we can observe that $\iota_Y$ induces an involution $\iota_Y'$ on $Y_6^{(1)'}$ and the two lifts of this involution on $S_6^{(1)}$ are $\sigma_2$ and $\sigma_3$.}
\end{rem}

\subsubsection{The surfaces $Y_9^{(2)}$ and $X_9^{(2)}$}\label{subsect: Y92 and X92}

Let $Y_9^{(2)}$ be a K3 surface such that $NS(Y_9^{(2)})\simeq L_9^{(2)}$.  With the same notation as Proposition \ref{prop: Y_n with NS=Ln2}, we have that $\varphi_{|c/2|}:Y_9^{(2)}\ra\mathbb{P}^2$ is the double cover branched along the conic $u:=\varphi_{|c/2|}\left(\left(c-\sum_{i=1}^8r_i\right)/2\right)$ and the quartic $q:=\varphi_{|c/2|}\left(c-\left(\sum_{i=1}^8r_i\right)/2\right)$. We denote by $\iota_Y$ the cover involution on $Y_9^{(2)}$. 
By Proposition \ref{prop: Y_n with NS=Ln2}
$Y_9^{(2)}$admits a double cover $X_9^{(2)}\ra Y_9^{(2)}$ which satisfies the Condition \ref{condition}. Let us now consider the curve $C\subset Y_9^{(2)}$ which is the curve of genus 5 contained in the branch locus of the double cover $X_9^{(2)}\ra Y_9^{(2)}$. Since the class of $C$ is $c$ and the double cover $Y_9^{(2)}\ra \mathbb{P}^2$ is defined by the linear system of the divisor $c/2$, the image of $C$ in $\mathbb{P}^2$, is the conic $v:=\varphi_{|c/2|}(C)\subset\mathbb{P}^2$. We observe that $u\neq v$, and indeed $u$ is in the branch locus of the double cover $\varphi_{|c/2|}:Y_9^{(2)}\ra \mathbb{P}^2$ and $v$ is not.

This situation is very similar to the one analyzed in Section \ref{subsec: Y6 and X6} (indeed only the degree of the curves in $\mathbb{P}^2$ changes) and so one obtains similar results: 

\begin{proposition}
Let $S_9^{(2)}$ be the bidouble cover of $\mathbb{P}^2$ branched on $u\cup v\cup q$ with multiplicity 2 on each component. Then the surface $X_9^{(2)}$ is birational to $S_9^{(2)}$.
\end{proposition}

By construction, we have the following commutative diagram:
$$\xymatrix{&S_9^{(2)}\ar[dr]\ar[d]\ar[dl]&\\
Y_9^{(2)'}\ar[dr]&W_9^{(2)'}\ar[d]&V_9^{(2)'}\ar[dl]\\
&\mathbb{P}^2&}$$
where $W_9^{(2)'}$ is the double cover of $\mathbb{P}^2$ branched along $v\cup q$ and $V_9^{(2)'}$ is the double cover of $\mathbb{P}^2$ branched along $u\cup v$. The minimal model $W_9^{(2)}$ of $W_9^{(2)'}$ is a K3 surface (since $W_9^{(2)'}$ is a double cover of $\mathbb{P}^2$ branched on a singular sextic) and the minimal model $V_9^{(2)}$ of $V_9^{(2)'}$ is a rational surface. In particular $h^{2,0}(W_9^{(2)})=1$ and $h^{2,0}(V_9^{(2)})=0$.
So also in this case the Hodge decomposition of $H^2(X_9^{(2)},\Z)$ splits in eigenspaces and in particular one finds two sub-Hodge structures of K3 type, both of them induced (by pull back) by the Hodge structure of the middle cohomology of the K3 surfaces $Y_9^{(2)}$ and $W_9^{(2)}$.

\begin{rem}{\rm As we already observed the construction of the surfaces $X_9^{(2)}$ and $X_{6}^{(1)}$ as bidouble cover of $\mathbb{P}^2$ is strongly related to the presence of the involution $\iota_Y$ on $Y_9^{(2)}$ and on $Y_{6}^{(1)}$ and thus to the fact that the lattice $NS(Y_9^{(2)})$ and $NS(Y_{6}^{(1)})$ are two elementary. These two lattices are the unique two elementary lattices which appears in the list $\mathcal{L}$.}\end{rem}

\subsection{The K3 surfaces $Y_n$ admit a lot of different double covers}
If $Y_{6}^{(1)}$ is a K3 surface with $NS(Y_{6}^{(1)})\simeq L_6^{(1)}$, then $Y_{6}^{(1)}$ admits the $2:1$ map $\varphi_{|C|}:Y_{6}^{(1)}\ra\mathbb{P}^2$, described in Proposition \ref{prop: Y_n with NS=Ln1}. In Section \ref{subsec: Y6 and X6} we described the geometry of the smooth double cover $X_{6}^{(1)}$ of $Y_{6}^{(1)}$ branched along the 2-divisible set $\{C, R_i\}$. We observe that the surface $Y_{6}^{(1)}$ admits also another smooth double cover, which is of the type described in Section \ref{subsect: NS(Y)=U(2)+D4}, i.e. branched along a smooth curve of genus 1 and 4 rational curves. Indeed, if one considers the pencil of lines in $\mathbb{P}^2$ through the point $\varphi_{|C|}(R_5)$ (which is a singular point of the branch locus of the double cover $\varphi_{|C|}:Y_{6}^{(1)}\ra \mathbb{P}^2$), this induces a pencil of curve of genus 1 on $Y_{6}^{(1)}$, whose class is $E\simeq C-R_5$. So $\varphi_{|E|}:Y_{6}^{(1)}\ra\mathbb{P}^1$ is a genus 1 fibration with a bisection, whose class is $R_5$. This fibration has a unique reducible fiber, which is of type $I_0^*$ and corresponds to the line $l\subset \mathbb{P}^2$ (i.e. with the unique line of the pencil of lines in $\mathbb{P}^2$ which is contained in the branch locus of the cover $\varphi_{|C|}:Y_{6}^{(1)}\ra\mathbb{P}^2$). The components of the fiber of type $I_0^*$ are the curves $R_i$, $i=1,\ldots, 4$, and the curve $D$, counted with multiplicity 2.

In the linear system $|E|$ there exists a smooth curve of genus 1, and we denote it by $E$. Let us now consider the set of disjoint curves $\{E, R_1,R_2,R_3,R_4\}$. This is a 2-divisible set on $Y_{6}^{(1)}$, indeed $E+R_1+R_2+R_3+R_4\simeq C-R_5+R_1+R_2+R_3+R_4\simeq (C+R_5+R_1+R_2+R_3+R_4)-2R_5\simeq 2(L-R_5)\in NS(Y_{6}^{(1)})$.
So one can construct a double cover of $Y_{6}^{(1)}$ (which is of course totally different from $X_{6}^{(1)}$) which is branched along $E$ and $R_i$, $i=1,2,3,4$. This is indeed what is done in Section \ref{subsect: NS(Y)=U(2)+D4}.

Now let us consider the K3 surfaces $Y_n$ such that $NS(Y_n)\simeq L_{n}^{(1)}$ for $n>6$. All these K3 surfaces admit a model as double cover of $\mathbb{P}^2$ branched along the union of a line and a quintic, by Proposition \ref{prop: projection of Y_n on P2}, which is associated to the linear system $|C-\sum_{i=1}^{n-6} R_i|$ (with the same notation of the proof of Proposition \ref{prop: projection of Y_n on P2}). So there is a smooth irreducible curve of genus 2 in the linear system $|C-\sum_{i=1}^{n-6} R_i|$ and we denote this curve by $M$. The set $\{M, R_{n-5},\ldots, R_{n-1}\}$ is a 2-divisible set of disjoint smooth curves, so $Y_n$ admits not only the double cover branched along $\{C, R_1,\ldots, R_{n-1}\}$ but also a double cover branched along $\{M, R_{n-5},\ldots, R_{n-1}\}$. This latter cover is of the type described in Section \ref{subsec: Y6 and X6}, i.e. it is branched along a smooth curve of genus 2 and 5 disjoint rational curves. Moreover, by the same argument as before (i.e. considering the pencil of lines through a singular point of the branch locus $\varphi_{|C-\sum_{i=1}^{n-6} R_i|}:Y_n\ra\mathbb{P}^2$) one proves that $Y_n$ admits also a double cover branched along a smooth curve of genus 1 and 4 rational curves, i.e. a double cover as in Section \ref{subsect: NS(Y)=U(2)+D4}.

Similarly, for all the admissible $n$, the K3 surfaces $Y_{n}^{(2)}$ such that $NS(Y_n^{(2)})\simeq L_n^{(2)}$ admit a double cover branched along the union of a smooth genus 1 curve and smooth 4 rational curves (as the double cover descried in Section \ref{subsect: NS(Y)=U(2)+D4}) and a double cover branched along a curve of genus 5 and 8 disjoint rational curves (as the one described in Section \ref{subsect:  Y92 and X92}). Indeed, by the Proposition \ref{prop: projection of Y_n(2) on P2}, the surfaces $Y_{n}^{(2)}$ admit a model as double covers of $\mathbb{P}^2$ branched along the union of a conic $u$ and a quartic $q$. Let us denote by $\{R_{n-8},\ldots R_{n-1}\}$ the curves contracted by this double covers and by $M$ a curve which is the pull back to $Y_n^{(2)}$ of a generic conic in $\mathbb{P}^2$. Then the set $\{M, R_{n-8},\ldots R_{n-1}\}$ is a 2-divisible set of smooth disjoint curves, so there exists a smooth double cover of $Y_n^{(2)}$ branched along this set. 

Let us now consider the pencil of conics in $\mathbb{P}^2$ passing through 4 of the eight points $u\cap q$ (we can assume that these points are the contraction of the curves $R_{n-8}, \ldots, R_{n-5}$). This pencil of conics induces a pencil of genus 1 curve on $Y_n^{(2)}$ and thus a genus 1 fibration $Y_n^{(2)}\ra\mathbb{P}^1$. We denote by $E$ a smooth fiber of this fibration. The set $\{E, R_{n-4}, \ldots R_{n-1}\}$ is a 2-divisible set of smooth curves on $Y_n^{(2)}$ and so $Y_n^{(2)}$ admits a smooth double cover branched along a curve of genus 1 and 4 rational curves.

More in general, if one proves that the linear system $|C-\sum_{i=1}^k R_i|$ contains a smooth irreducible curve $M$, then it has genus $n-4-k$ and $\{ M, R_{k+1},\ldots R_{n-1}\}$ is a 2-divisible set which allows one to construct a double cover of $Y_n^{(r)}$ branched along a curve of genus  $n-4-k$  and $n-k$ rational curves. 

\section{Remarks on the cases $g(C)>1$, $h\neq -1$}\label{sec: double cover h not -1}
In this section we would like to describe some generalizations of the previous results to the case $h\neq -1$. In particular we will investigate the existence of K3 surfaces admitting a double cover as in Condition \ref{cond 2} providing results both in the case of the biggest possible value of $n$ and in the case of the lowest possible value of $n$. Then we discuss an example, which is a generalization of the construction described in Section \ref{subsec: bidouble covers P2}. 
\begin{condition}\label{cond 2}
Let $S_n^{(h)}$ be a K3 surface with $\rho(S_n^{(h)})=n$ and let $Z_n^{(h)}\ra S_n^{(h)}$ be a smooth double cover branched on $n$ curves such that $n-1$ of them, denoted by $R_i$ $i=1,\ldots, n-1$, are rational and the the remaining curve, denoted by $C$, has genus $g(C)=n+4h$. We assume $h\geq -3$ and $g(C)\geq 2$.
\end{condition}

\begin{proposition}\label{prop: existence h not -1}
If $n=17$ there exist K3 surfaces as in Condition \ref{cond 2} if and only if $h$ is even. In this case $S_{17}^{(h)}\simeq Km(A^{(h)})$ where $A^{(h)}$ is an Abelian surface with a polarization of degree $4+h$.

If $n=16$ there exists a K3 surface as in Condition \ref{cond 2} for every $h\geq -3$.

If $h\geq 1$ there exists a K3 surface as in Condition \ref{cond 2} with $n=1$.
\end{proposition}
\proof Let us assume that $n=17$ and that $S_{17}^{(h)}$ satisfies Condition \ref{cond 2}. Then, by construction, the classes of the curves $C$ and $R_i$, $i=1,\ldots,16$ span, over $\Q$, the N\'eron--Severi group of the surface $S_{17}^{(h)}$. By \cite{NikKummer}, if there are 16 disjoint rational curves on a K3 surface, then the minimal primitive sublattice of the N\'eron--Severi group which contains these curves is the Kummer lattice $K$. So the lattice $K$ is primitively embedded in $NS(S_{17}^{(h)})$, and hence $S_{17}^{(h)}$ is a Kummer surface $Km(A^{(h)})$ for a certain Abelian surface $A^{(h)}$, cf.\ \cite{NikKummer}. Moreover, the class $C$ is orthogonal to $K$ and the class $(C+\sum_{i=1}^{16} R_i)/2$ is contained in $NS(S_{17}^{(h)})$. Since $(\sum_{i=1}^{16} R_i)/2\in K\subset NS(S_{17}^{(h)})$, it follows that $C/2\in NS(S_{17}^{(h)})$. So $S_{17}^{(h)}$ has a polarization, orthogonal to $K$ which has degree $(C/2)^2=(g(C)-1)/2=8+2h$. Since the orthogonal to $K$ in $NS(Km(A^{(h)}))$ is isometric to $NS(A^{(h)})(2)$ and $NS(A^{(h)})$ is an even lattice, it follows that $8+2h\equiv 0\mod 4$, i.e. $h$ is even. Moreover, since on $Km(A^{(h)})$ we have a polarization ($C/2$) with degree $8+2h$, on $A^{(h)}$ there is a polarization of degree $4+h$.

Let us now consider an Abelian surface $A^{(h)}$ such that $NS(A^{(h)})\simeq 4+h$ and let us consider the Kummer surface $Km(A^{(h)})$. This is a K3 surface on which there are 16 rational curves, denoted by $R_i$, $i=1,\ldots, 16$ arising from the desingularization of the quotient $A^{(h)}/\iota_A$. Moreover there is a polarization, say $H$, induced by the polarization on $A^{(h)}$ and with degree $8+2h$. The divisor $H$ is pseudoample and without fixed components, thus a general member of the linear system $|2H|$ is smooth and will be denoted by $C$. Now it is clear that $C+\sum_{i=1}^{16}R_i$ is a 2-divisible set, which allows one to construct a double cover of $S_{17}^{(h)}\simeq Km(A^{(h)})$ as in Condition \ref{cond 2}.

The case $n=16$ is similar: one first observe that the curves $R_i$, $i=1,\ldots, 15$ spans the lattice $M_{(\Z/2\Z)^4}$ and then one observes that in this case $C^2=2g(C)-2=30+8h=2(15+4h)$. Denoted by  $d:=15+4h$ one obtains that $d\equiv 3\mod 4$. Then \cite[Theorem 8.3]{GS} implies that there exists a K3 surface $S_{16}^{(h)}$ whose N\'eron--Severi group is an overlattice of index 2 of $\langle 2d\rangle\oplus M_{(\Z/2\Z)^4}$ obtained by adding exactly the vector $(C+\sum_{i=1}^{15} R_i)/2$ to $\langle 2d\rangle\oplus M_{(\Z/2\Z)^4}$. But this is the vector which is associated to the required double cover. 
Viceversa, let $S_{16}^{(h)}$ be a K3 surface whose N\'eron--Severi group is an overlattice of index 2 of $\langle 2d\rangle\oplus M_{(\Z/2\Z)^4}$ obtained adding $(C+\sum_{i=1}^{15} R_i)/2$ to $\langle 2d\rangle\oplus M_{(\Z/2\Z)^4}$, where $C$ is the class spanning  $\langle 2d\rangle$ and the $R_i$'s are a basis of the roots lattice of $M_{(\Z/2\Z)^4}$. Then, by \cite[Theorem 5.2]{G}, the $R_i$'s represent smooth irreducible curves and $C$ is a polarization, whose linear system contains a smooth irreducible curve. Since $\{C,R_i\}$ is an even set, by the properties of the N\'eron--Severi group, we conclude that the K3 surface $S_{16}^{(h)}$ is as in Condition \ref{cond 2}.

The case $n=1$ and $h\geq 1$ was already discussed in Section \ref{sebsect:  existence}.\endproof

\subsection{The cases $h=0$ and $h=-2$, an example}
In Sections \ref{subsec: Y6 and X6} and \ref{subsect:  Y92 and X92} we considered two K3 surfaces which are double covers of $\mathbb{P}^2$ branched along a sextic $B$ which is reducible and has exactly two components, in one case a line and a quintic in the other case a conic and a quartic. Both these K3 surfaces admit a smooth double cover and it turns out that the smooth double cover is birational to a bidouble cover of $\mathbb{P}^2$ branched along the union of three curves, two of which are the components of the reducible sextic $B$. 
It is now natural to consider the unique remaining case of sextic which splits in two components, the case of a sextic $B$ which is the union of two cubics $b_1$ and $b_2$. In this section we consider the K3 surface $Y$ which is a double cover of $\mathbb{P}^2$ branched on $B=b_1\cup b_2$. We show that $Y$ admits two double covers ($X$ and $Z$) as in Condition \ref{cond 2} (for different values of $h$) and both these double covers are bidouble covers of $\mathbb{P}^2$ branched along $b_1$, $b_2$ and a third curve.

Let us consider a K3 surface $Y$ whose N\'eron--Severi group is the lattice generated over $\Q$ by the classes $c$, $r_i$, $i=1,\ldots 9$ with $c^2=2$, $r_i^2=-2$, $cr_i=r_ir_j=0$, $i\neq j$ and which is generated over $\Z$ by the same classes and by $(c+\sum_{i=1}^9 r_i)/2$. Then $NS(Y)$ is a 2-elementary lattice with rank 10, length 8 and whose discriminant form takes values in $\Z$. The divisor $c$ can be chosen to be pseudo ample and $\varphi_{|c|}:Y\ra \mathbb{P}^2$ is a generically $2:1$ cover branched along the union of two cubics $b_1$ and $b_2$. The class of each component of the ramification curve in $Y$ is $(3c-\sum_{i=1}^9r_i)/2$ and the intersection points $b_1\cap b_2$ are $\varphi_{|c|}(r_i)$.

In the linear system $|c|$ there is a a smooth irreducible curve $C_2$ of genus 2, and we denote by $l$ the line $\varphi_{|c|}(C_2)$. Similarly, in the linear system $|3c|$ there is a smooth irreducible curve $C_{10}$ of genus 10 and we denote by $b_3$ the cubic $\varphi_{|c|}(C_{10})$.

By construction the set $\{C_2, R_1,\ldots R_9\}$ is 2-divisible on $Y$ and the smooth surface $X$ which is the double cover of $Y$ branched along $C_2\bigcup \cup_{i=1}^9 R_i$ is an example of smooth double cover of a K3 surface as in Condition \ref{cond 2} with $n=10$ and $h=-2$.

Similarly, the set $\{C_{10}, R_1,\ldots R_9\}$ is 2-divisible on $Y$ and the smooth surface $Z$ which is the double cover of $Y$ branched along $C_{10}\bigcup \cup_{i=1}^9 R_i$ is an example of smooth double cover of a K3 surface as in Condition \ref{cond 2}  with $n=10$ and $h=0$.

\begin{proposition}
Let $U_2$ be the bidouble cover of $\mathbb{P}^2$ branched on $ b_1\cup b_2\cup l$ with multiplicity 2 on each component. Then the surface $X$ is birational to $U_2$.

Let $U_{10}$ be the bidouble cover of $\mathbb{P}^2$ branched on $b_1\cup b_2\cup b_3$ with multiplicity 2 on each component. Then the surface $Z$ is birational to $U_{10}$.
\end{proposition}

In particular, we observe that there exists the following diagram:
$$
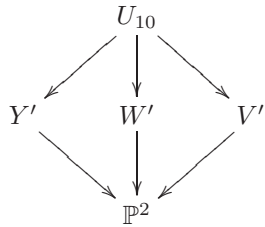
\xymatrix{&U_{10}\ar[dr]\ar[d]\ar[dl]&\\
Y'\ar[dr]&W'\ar[d]&V'\ar[dl]\\
&\mathbb{P}^2&}$$
where $W'$ is the double cover of $\mathbb{P}^2$ branched along $b_2\cup b_3$ and $V'$ is the double cover of $\mathbb{P}^2$ branched along $b_1\cup b_3$. So the minimal models of the surfaces $Y'$, $W'$, $V'$ are all K3 surfaces, denoted by $Y$, $W$ and $V$. We have $h^{2,0}(U_{10})=h^{2,0}(Z)=3$ and the Hodge structure of $U_{10}$ (and thus of $Z$) splits in three sub-Hodge structures of K3 type, each of them induced by the Hodge structure of a K3 surface (i.e. of $Y$, $W$ and $V$).

\end{document}